\documentclass[11pt]{amsart}
\usepackage{tikz}
\usepackage{amsmath, amssymb, amsthm, mathtools}
\usepackage{mathrsfs}
\usepackage[T1]{fontenc}
\usepackage{microtype}
\AddToHook{env/thebibliography/begin}{\renewcommand{\etalchar}[1]{}}
\usepackage{enumitem}
\usepackage{array, booktabs}
\usepackage[colorlinks=true, linkcolor=blue!50!black, citecolor=green!40!black,
            urlcolor=blue!60!black]{hyperref}
\usepackage[capitalize, nameinlink]{cleveref}
\theoremstyle{plain}
\newtheorem{theorem}{Theorem}[section]
\newtheorem{proposition}[theorem]{Proposition}
\newtheorem{lemma}[theorem]{Lemma}
\newtheorem{corollary}[theorem]{Corollary}
\theoremstyle{definition}
\newtheorem{remark}[theorem]{Remark}
\newtheorem{example}[theorem]{Example}

\newcommand{\Z}{\mathbb{Z}}
\newcommand{\M}{\mathsf{M}}                 % matching complex
\newcommand{\Mp}{\mathsf{M}_p}              % perfect matching complex
\newcommand{\Ind}{\mathsf{Ind}}             % independence complex
\newcommand{\AP}{\mathrm{AP}}               % antiprism graph
\newcommand{\GP}{\mathrm{GP}}               % generalized Petersen graph

\tikzset{
  gpvertex/.style = {circle, fill=white, draw=black, line width=.5pt, inner sep=0pt,
                     minimum size=1.8mm},
  gpblack/.style  = {gpvertex, fill=black},
  gpgrey/.style   = {gpvertex, fill=black!30},
  gpedge/.style   = {line width=.5pt, black!60},
  gpinner/.style  = {gpedge, bend right=35},
  gptoinf/.style  = {gpedge, -stealth},
}
\makeatletter
\newcommand\gp@list[2]{\gdef#1{,}\foreach\gp@x in {#2}{\xdef#1{#1\gp@x,}}}
\newcommand\gp@mem[3]{%
  \edef\gp@tmp{\noexpand\in@{,#1,}{#2}}\gp@tmp\ifin@\def#3{1}\else\def#3{0}\fi}
\newcommand\gpgraph[4][]{%
  \gp@list\gp@U{#3}\gp@list\gp@V{#4}\def\gp@opt{,#1,}%
  \gdef\gp@first{}%
  \foreach\gp@x in {#2}{\ifx\gp@first\@empty\xdef\gp@first{\gp@x}\fi\xdef\gp@last{\gp@x}}%
  \foreach\gp@i in {#2}{%
    \pgfmathtruncatemacro\gp@j{\gp@i+1}\pgfmathtruncatemacro\gp@k{\gp@i+2}%
    \gp@mem\gp@i\gp@U\gp@a\gp@mem\gp@j\gp@U\gp@b\gp@mem\gp@i\gp@V\gp@c\gp@mem\gp@k\gp@V\gp@d
    \ifnum\gp@a\gp@b=11 \draw[gpedge] (\gp@i,1) -- (\gp@j,1);\fi
    \ifnum\gp@a\gp@c=11 \draw[gpedge] (\gp@i,1) -- (\gp@i,0);\fi
    \ifnum\gp@c\gp@d=11 \draw[gpinner] (\gp@i,0) to (\gp@k,0);\fi}%
  \gp@mem{right}\gp@opt\gp@a
  \ifnum\gp@a=1
    \gp@mem\gp@last\gp@U\gp@b\ifnum\gp@b=1 \draw[gptoinf] (\gp@last,1) -- ++(0.8,0);\fi
    \pgfmathtruncatemacro\gp@j{\gp@last-1}%
    \foreach\gp@x in {\gp@j,\gp@last}{%
      \gp@mem\gp@x\gp@V\gp@b\ifnum\gp@b=1 \draw[gptoinf] (\gp@x,0) -- ++(-35:5mm);\fi}%
  \fi
  \gp@mem{left}\gp@opt\gp@a
  \ifnum\gp@a=1
    \gp@mem\gp@first\gp@U\gp@b\ifnum\gp@b=1 \draw[gptoinf] (\gp@first,1) -- ++(-0.8,0);\fi
    \pgfmathtruncatemacro\gp@j{\gp@first+1}%
    \foreach\gp@x in {\gp@first,\gp@j}{%
      \gp@mem\gp@x\gp@V\gp@b\ifnum\gp@b=1 \draw[gptoinf] (\gp@x,0) -- ++(215:5mm);\fi}%
  \fi
  \foreach\gp@x in {#3}{\node[gpvertex] at (\gp@x,1) {};}%
  \foreach\gp@x in {#4}{\node[gpvertex] at (\gp@x,0) {};}}
\newcommand\gpcols[2][-0.75]{%
  \foreach\gp@x in {#2}{\node[font=\scriptsize,text=black!55] at (\gp@x,#1) {$\gp@x$};}}
\tikzset{
  gpgone/.style       = {line width=.4pt, black!20},
  gpkeep/.style       = {line width=.9pt, black},
  gpgonevertex/.style = {circle, fill=black!20, inner sep=0pt, minimum size=1.2mm},
}
\newcommand\gpcircle[3]{%
  \gp@list\gp@U{#2}\gp@list\gp@V{#3}%
  \pgfmathtruncatemacro\gp@n{#1-1}%
  \foreach\gp@i in {0,...,\gp@n}{%
    \coordinate (u\gp@i) at ({90-360*\gp@i/#1}:1.3);%
    \coordinate (v\gp@i) at ({90-360*\gp@i/#1}:0.65);}%
  \foreach\gp@i in {0,...,\gp@n}{%
    \pgfmathtruncatemacro\gp@j{mod(\gp@i+1,#1)}\pgfmathtruncatemacro\gp@k{mod(\gp@i+2,#1)}%
    \draw[gpgone] (u\gp@i)--(u\gp@j) (u\gp@i)--(v\gp@i) (v\gp@i)--(v\gp@k);}%
  \foreach\gp@i in {0,...,\gp@n}{%
    \pgfmathtruncatemacro\gp@j{mod(\gp@i+1,#1)}\pgfmathtruncatemacro\gp@k{mod(\gp@i+2,#1)}%
    \gp@mem\gp@i\gp@U\gp@a\gp@mem\gp@j\gp@U\gp@b\gp@mem\gp@i\gp@V\gp@c\gp@mem\gp@k\gp@V\gp@d
    \ifnum\gp@a\gp@b=11 \draw[gpkeep] (u\gp@i)--(u\gp@j);\fi
    \ifnum\gp@a\gp@c=11 \draw[gpkeep] (u\gp@i)--(v\gp@i);\fi
    \ifnum\gp@c\gp@d=11 \draw[gpkeep] (v\gp@i)--(v\gp@k);\fi}%
  \foreach\gp@i in {0,...,\gp@n}{%
    \gp@mem\gp@i\gp@U\gp@a\gp@mem\gp@i\gp@V\gp@c
    \ifnum\gp@a=1
      \node[gpvertex] at (u\gp@i) {};
      \node[font=\scriptsize] at ({90-360*\gp@i/#1}:1.62) {$u_{\gp@i}$};
    \else\node[gpgonevertex] at (u\gp@i) {};\fi
    \ifnum\gp@c=1
      \node[gpvertex] at (v\gp@i) {};
      \node[font=\scriptsize] at ({90-360*\gp@i/#1}:0.33) {$v_{\gp@i}$};
    \else\node[gpgonevertex] at (v\gp@i) {};\fi}}
\makeatother
\newcommand\petersenpm[2]{%
  \begin{tikzpicture}[scale=0.55,baseline=(current bounding box.center)]
    \foreach \i in {0,...,4}{%
      \coordinate (u\i) at ({90-72*\i}:1.4);
      \coordinate (v\i) at ({90-72*\i}:0.7);}
    \foreach \i in {0,...,4}{%
      \pgfmathtruncatemacro{\ia}{mod(\i+1,5)}\pgfmathtruncatemacro{\ib}{mod(\i+2,5)}%
      \draw[black!30,line width=0.6pt] (u\i)--(u\ia) (u\i)--(v\i) (v\i)--(v\ib);}
    \foreach \a/\b in {#1}{\draw[line width=1.6pt] (\a)--(\b);}
    \foreach \i in {0,...,4}{\fill (u\i) circle (1.6pt); \fill (v\i) circle (1.6pt);}
    \node[font=\small] at (0,-1.95) {#2};
  \end{tikzpicture}}

\title[Antiprism graphs and generalized Petersen graphs]{Perfect matching complexes\\ of antiprism
  graphs and independence\\ complexes of $\GP(n,2)$}
\author{Youngsoo Kim}
\address{Department of Mathematics, Tuskegee University, Tuskegee, AL 36088, USA}
\email{ykim1@tuskegee.edu}

\begin{document}

\begin{abstract}
We show that the perfect matching complex of the antiprism graph, the square of the $2m$-cycle,
is homotopy equivalent to the suspension of the independence complex of the cube of the same
cycle, wedged with four $(m/2-1)$-spheres when $m$ is even.
With Adamaszek's recursion and the initial values we compute, all these complexes are explicit
wedges of spheres.
We also show that the independence complex of the generalized Petersen graph $\GP(n,2)$ is a
sphere or a wedge of four spheres, and that it is periodic in $n$ with period $5$ up to a
threefold suspension.
\end{abstract}

\maketitle

\section{Introduction}\label{sec:intro}

Simplicial complexes built from graphs have been studied extensively. See the survey of Wachs
\cite{Wachs} and the book of Jonsson \cite{Jonsson}.
Two of the most studied are the following.
A set of vertices of a graph $G$ is \emph{independent} if no two of them are adjacent, and the
\emph{independence complex} $\Ind(G)$ is the simplicial complex whose faces are the independent
sets of $G$.
A \emph{matching} of $G$ is a set of edges, no two of which share an endpoint, and the
\emph{matching complex} $\M(G)$ is the simplicial complex whose faces are the matchings of $G$.
The \emph{line graph} of a graph has the edges of the graph as its vertices, and two of them
are adjacent if they share an endpoint.
So the matching complex of a graph is the independence complex of its line graph.
A simplicial complex is a \emph{flag complex} (or \emph{clique complex}) if a set $\sigma$ is a
face if and only if $\{a,b\}$ is a face for all $a,b\in\sigma$.
Independence complexes are flag complexes.
Kozlov \cite{Kozlov99} determined the homotopy types of the independence complexes of paths
and cycles.
For independence complexes there are now effective general tools, such as the fold lemma
\cite{Engstrom09} and the cofibre sequence of a vertex.
With these tools the homotopy types have been determined for powers of paths
\cite{Engstrom09} and for powers of cycles \cite{Adamaszek}.

The \emph{perfect matching complex} $\Mp(G)$ of a graph $G$ is the subcomplex of $\M(G)$
generated by the perfect matchings, the matchings that cover every vertex.
It was introduced by Bayer, Jeli\'c Milutinovi\'c and Vega \cite{BJVhoney}, and its homotopy
type is known for a few families of graphs: some honeycomb graphs \cite{BJVhoney}, polygonal
line tilings \cite{PLTpm}, circular and M\"obius ladders \cite{AB}, and the $3\times 2n$ grid
graphs \cite{CCGKS}.
The tools above for independence complexes do not apply directly to perfect matching complexes,
because a perfect matching complex is in general not a flag complex.
The works above get around this in different ways.
Bayer, Jeli\'c Milutinovi\'c and Vega \cite{BJVhoney} use a bijection between the perfect
matchings of honeycomb graphs and plane partitions, together with discrete Morse theory.
Chandrakar and Singh \cite{PLTpm}, and Agarwal and Basak \cite{AB}, use the minimal matchings
that are not contained in a perfect matching to relate the perfect matching complex to an
independence complex.
Chandrakar et al.\ \cite{CCGKS} cover the perfect matching complex by the stars of the edges at
a corner vertex and argue by induction.

In this paper we study two families of graphs: the antiprism graphs $\AP_m$ and the
generalized Petersen graphs $\GP(n,2)$.
For $\AP_m$ the independence complex is known \cite{Adamaszek}, and we determine the perfect
matching complex.
For $\GP(n,2)$ we determine the independence complex, which is periodic in $n$, and we find
that the perfect matching complexes show no such pattern.

The antiprism graph $\AP_m$ is the square $C_{2m}^2$ of the cycle of length $2m$
(\cref{fig:graphs}).
The graph $\AP_m$ is also obtained by joining the two ends of a line tiling of triangles, whose
perfect matching complex is determined in \cite[Thm.~5.3]{PLTpm}.
Agarwal and Basak ask for the independence and perfect matching complexes of other families of
graphs related to the circular and M\"obius ladders, in particular of $4$-regular ones
\cite[Problem~27]{AB}, and $\AP_m$ is such a graph.
The generalized Petersen graphs $\GP(n,2)$ (\cref{fig:graphs}) are related to the ladders in
another way: like them, they are $3$-regular, and $\GP(n,1)$ is the circular ladder.
They include the Petersen graph $\GP(5,2)$ and the graph $\GP(10,2)$ of the dodecahedron.
The independence complex of $\GP(n,1)$ is determined in \cite{AB}, and that of the Petersen
graph by Barmak \cite{Barmak13}.

Our main results are the following.
First, the perfect matching complex of the square of a cycle is determined by the independence
complex of the \emph{cube} of the same cycle: for $m\ge3$ we show (\cref{thm:main})
\[
  \Mp(\AP_m)\simeq
  \begin{cases}
    \Sigma\Ind(C_{2m}^3) & \text{if $m$ is odd},\\[2pt]
    \bigvee_4 S^{m/2-1}\ \vee\ \Sigma\Ind(C_{2m}^3) & \text{if $m$ is even}.
  \end{cases}
\]
So the perfect matching complexes of these graphs are governed by the independence complexes of
powers of cycles, and the results of Adamaszek \cite{Adamaszek} on the latter apply to them.
Second, the independence complexes of the graphs $\GP(n,2)$ are periodic in $n$, with period
$5$ and a threefold suspension in each period: for $n\ge5$, written as $n=5q+r$ with
$0\le r\le4$, we show (\cref{thm:gp})
\[
  \Ind(\GP(n,2))\simeq
  \begin{cases}
    \bigvee_4 S^{3q-1} & \text{if $r=0$},\\[2pt]
    S^{3q-2+r} & \text{if $1\le r\le4$}.
  \end{cases}
\]

In the proof of \cref{thm:main} we get around the fact that $\Mp(\AP_m)$ is not a flag complex
differently from the works above: we split it into two flag complexes.
A perfect matching of $\AP_m$ divides the cycle into intervals of two or four consecutive
vertices, and the first indices of these intervals all have the same parity.
So $\Mp(\AP_m)$ is the union of two subcomplexes $K_0$ and $K_1$, one for each parity.
These are flag complexes, so the tools for independence complexes apply to them.
Their intersection is $\Ind(C_{2m}^3)$, and we glue them along it.

Adamaszek gave a recursion for $\Ind(C_n^r)$ \cite{Adamaszek}, but not its initial values.
For $r=3$ we compute them with SageMath \cite{Sage} in \cref{sec:explicit}, and with the
recursion they show that
$\Ind(C_n^3)$ is a wedge of spheres for every $n\ge6$.
By \cref{thm:main}, so is $\Mp(\AP_m)$ for every $m\ge3$, and the number of spheres of each
dimension in both is determined (\cref{cor:wedge}).

The proof of \cref{thm:gp} needs only the standard tools for independence complexes.
Three cofibre sequences cut $\GP(n,2)$ open into strips, and removing five columns of a strip
suspends its independence complex three times (\cref{prop:block}).
In \cref{sec:gp-pm} we present the integral homology of $\Mp(\GP(n,2))$ for $n\le14$, computed
with SageMath \cite{Sage}, and it shows no evident pattern (\cref{tab:gp-pm}).
For the dodecahedron $\GP(10,2)$ it has torsion.

\subsection{Notation}\label{sec:notation}

All graphs in this paper are finite and simple.
We write $uv$ (or $vu$) for the edge between vertices $u$ and $v$, $N(v)$ for the set of
neighbours of a vertex $v$, and $N[v]=N(v)\cup\{v\}$.
If $G$ has a perfect matching, the facets of $\Mp(G)$ are its perfect matchings, and the faces
of $\Mp(G)$ are the matchings contained in a perfect matching.
A matching of $d+1$ edges is a $d$-simplex.
We reserve the words \emph{vertex} and \emph{edge} for graphs, so a $0$-simplex of $\Mp(G)$ is
called an edge.
We write $\bigvee_k X$ for the wedge of $k$ copies of a space $X$, and $\Sigma X=S^0*X$ for
the unreduced suspension.
We write $S^{-1}$ for the complex $\{\emptyset\}$, whose geometric realization is empty, so
that $\Sigma S^{-1}=S^0$. For example, the independence complex of the empty graph is $S^{-1}$.

\subsection{The graphs}\label{sec:graphs}

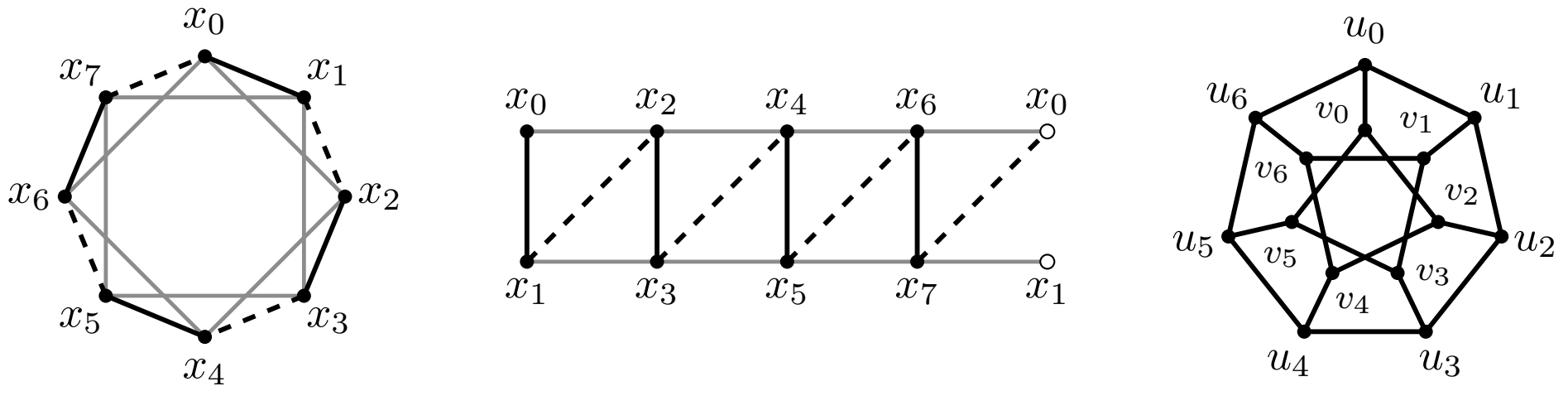
\begin{figure}[ht]
\centering
\begin{tikzpicture}[scale=0.8,baseline=(current bounding box.center)]
  \foreach \i in {0,...,7}{\coordinate (c\i) at ({90-45*\i}:1.4);}% vertex x_i
  % edges l_i = x_i x_{i+2} (grey)
  \foreach \i in {0,...,7}{%
    \pgfmathtruncatemacro{\ii}{mod(\i+2,8)}%
    \draw[black!45,line width=0.9pt] (c\i)--(c\ii);}
  % edges s_i = x_i x_{i+1}: s_{2j} solid, s_{2j+1} dashed
  \foreach \j in {0,...,3}{%
    \pgfmathtruncatemacro{\ia}{2*\j}\pgfmathtruncatemacro{\ib}{2*\j+1}%
    \pgfmathtruncatemacro{\ic}{mod(2*\j+2,8)}%
    \draw[line width=1.1pt] (c\ia)--(c\ib);
    \draw[line width=1.1pt,dashed] (c\ib)--(c\ic);}
  \foreach \i in {0,...,7}{%
    \fill (c\i) circle (2pt);
    \node[font=\small] at ({90-45*\i}:1.75) {$x_{\i}$};}
\end{tikzpicture}
\hfill
\begin{tikzpicture}[scale=0.8,baseline=(current bounding box.center)]
  \foreach \j in {0,...,4}{%
    \coordinate (u\j) at (1.3*\j,1.3);
    \coordinate (w\j) at (1.3*\j,0);}
  \foreach \j in {0,...,3}{%
    \pgfmathtruncatemacro{\jn}{\j+1}%
    \draw[black!45,line width=0.9pt] (u\j)--(u\jn);
    \draw[black!45,line width=0.9pt] (w\j)--(w\jn);
    \draw[line width=1.1pt] (u\j)--(w\j);
    \draw[line width=1.1pt,dashed] (w\j)--(u\jn);}
  \foreach \j in {0,...,3}{%
    \pgfmathtruncatemacro{\ev}{2*\j}\pgfmathtruncatemacro{\od}{2*\j+1}%
    \fill (u\j) circle (2pt); \fill (w\j) circle (2pt);
    \node[font=\small,above] at (u\j) {$x_{\ev}$};
    \node[font=\small,below] at (w\j) {$x_{\od}$};}
  % the last column repeats the first
  \filldraw[fill=white,draw=black] (u4) circle (2pt);
  \filldraw[fill=white,draw=black] (w4) circle (2pt);
  \node[font=\small,above] at (u4) {$x_0$};
  \node[font=\small,below] at (w4) {$x_1$};
\end{tikzpicture}
\hfill
\begin{tikzpicture}[scale=0.8,baseline=(current bounding box.center)]
  \foreach \i in {0,...,6}{%
    \coordinate (u\i) at ({90-360*\i/7}:1.4);% outer vertex u_i
    \coordinate (v\i) at ({90-360*\i/7}:0.75);}% inner vertex v_i
  \foreach \i in {0,...,6}{%
    \pgfmathtruncatemacro{\ia}{mod(\i+1,7)}\pgfmathtruncatemacro{\ib}{mod(\i+2,7)}%
    \draw[line width=1.1pt] (u\i)--(u\ia) (u\i)--(v\i) (v\i)--(v\ib);}
  \foreach \i in {0,...,6}{%
    \fill (u\i) circle (2pt); \fill (v\i) circle (2pt);
    \node[font=\small] at ({90-360*\i/7}:1.75) {$u_{\i}$};
    \node[font=\scriptsize] at ({90-360*\i/7+19}:0.98) {$v_{\i}$};}
\end{tikzpicture}
\caption{The graphs $\AP_4$ (left and middle) and $\GP(7,2)$ (right). In $\AP_4$ the short
edges $s_i$ are solid for even $i$ and dashed for odd $i$, and the long edges $\ell_i$ are
grey.}\label{fig:graphs}
\end{figure}

For $n\ge3$ and $r\ge1$, the \emph{$r$-th power} $C_n^r$ of the $n$-cycle has vertices $x_i$
($i\in\Z_n$), and $x_i$ and $x_j$ are adjacent if and only if
$j-i\equiv\pm s\pmod n$ for some $1\le s\le r$, that is, if their distance in the cycle is at
most $r$.
Similarly, for $n\ge1$ the $r$-th power $P_n^r$ of the path with $n$ vertices has vertices
$x_1,\dots,x_n$, and $x_i$ and $x_j$ are adjacent if and only if $1\le|j-i|\le r$.
For $m\ge3$, the \emph{antiprism graph} is the square of the $2m$-cycle,
\[
  \AP_m=C_{2m}^2,
\]
the graph of the antiprism over an $m$-gon (\cref{fig:graphs}).
It is $4$-regular, with $2m$ vertices and $4m$ edges.
For $i\in\Z_{2m}$ we write
\[
  s_i=x_ix_{i+1},\qquad \ell_i=x_ix_{i+2}
\]
for the edges of $\AP_m$, and we call the edges $s_i$ \emph{short} and the edges $\ell_i$
\emph{long}.

The antiprism graph $\AP_m$ is also obtained from the circular ladder (or prism graph) on $2m$
vertices by adding a diagonal to each of its $m$ quadrilaterals.
The circular and the M\"obius ladders \cite{AB} are the graphs of an annulus and of a M\"obius
band made of quadrilaterals. Similarly, for $n\ge5$ the graph $C_n^2$ is the graph of a band
made of $n$ triangles instead of quadrilaterals.
It is an annulus if $n$ is even and a M\"obius band if $n$ is odd.
In the odd case, $C_n^2$ has an odd number of vertices, so it has no perfect matching.

The circular ladder is also the first member of a second family.
For $n\ge3$ and $1\le k<n/2$, the \emph{generalized Petersen graph} $\GP(n,k)$
\cite{Watkins69} has outer
vertices $u_i$ and inner vertices $v_i$ ($i\in\Z_n$), and edges
\[
  u_iu_{i+1},\qquad u_iv_i,\qquad v_iv_{i+k}\qquad(i\in\Z_n)
\]
(\cref{fig:graphs}). It is $3$-regular, with $2n$ vertices and $3n$ edges, and $\GP(n,1)$ is the circular ladder
on $2n$ vertices.
The graph $\GP(5,2)$ is the Petersen graph and $\GP(10,2)$ is the graph of the dodecahedron.
If $kl\equiv\pm1\pmod n$, then $\GP(n,k)\cong\GP(n,l)$ \cite[Sect.~2]{Watkins69}.
For distinct $k$ and $l$ coprime to $n$, the converse also holds \cite[Theorem~10]{SteimleStaton}.

\section{The perfect matching complex of the antiprism graph}\label{sec:antiprism}

In this section we prove \cref{thm:main}, which reduces $\Mp(\AP_m)$ to $\Ind(C_{2m}^3)$.
Throughout the section, we assume $m\ge3$ and indices are in $\Z_{2m}$.

\begin{lemma}\label{lem:tilings}
Every perfect matching $P$ of $\AP_m$ is a disjoint union of sets $\{s_i\}$ and
$\{\ell_j,\ell_{j+1}\}$, and every edge $\ell_i\in P$ lies in exactly one of these sets,
$\{\ell_{i-1},\ell_i\}$ or $\{\ell_i,\ell_{i+1}\}$.
Their vertex sets $\{x_i,x_{i+1}\}$ and $\{x_j,\dots,x_{j+3}\}$ partition the vertex set of
$\AP_m$.
Conversely, every such partition determines a perfect matching.
\end{lemma}

\begin{proof}
If $\ell_i=x_ix_{i+2}\in P$, then the vertex $x_{i+1}$ is matched with $x_{i-1}$ or with $x_{i+3}$,
but not both.
So exactly one of the sets $\{\ell_{i-1},\ell_i\}$ and $\{\ell_i,\ell_{i+1}\}$ is contained in $P$.
Write it as $\{\ell_j,\ell_{j+1}\}$ with $j\in\{i-1,i\}$. It matches the four consecutive vertices
$x_j,\dots,x_{j+3}$ among themselves.
Thus every vertex lies either in the vertex set $\{x_i,x_{i+1}\}$ of an edge $s_i\in P$, or in the
vertex set $\{x_j,\dots,x_{j+3}\}$ of a set $\{\ell_j,\ell_{j+1}\}\subseteq P$.
These vertex sets are pairwise disjoint because $P$ is a matching, so they partition the vertex
set of $\AP_m$.
The converse is clear.
\end{proof}

The vertex sets in \cref{lem:tilings} have even size and cover the cycle $x_0,\dots,x_{2m-1}$
consecutively, so their first indices are congruent modulo~$2$.
The \emph{parity} of $P$ is their common residue $p\in\{0,1\}$.
So if $P$ has parity $p$ and contains $\ell_i$, 
then $P$ contains the pair $\ell_i, \ell_{i+1}$ if $i\equiv p$, 
and $\ell_{i-1}, \ell_i$ if $i\not\equiv p$.
For $p\in\{0,1\}$ let $K_p$ be the subcomplex of $\Mp(\AP_m)$ generated by the perfect matchings
of parity $p$, so that $\Mp(\AP_m)=K_0\cup K_1$.
The rotation $x_i\mapsto x_{i+1}$ of $\AP_m$ shifts all indices by one, so it exchanges the two
parities, and $K_0\cong K_1$.
We now name the vertex sets of the partition in \cref{lem:tilings} for parity $p$.
For $i\in\Z_{2m}$ let
\[
  V_p(s_i)=\{x_i,x_{i+1}\},\qquad
  V_p(\ell_i)=
  \begin{cases}
    \{x_i,\dots,x_{i+3}\} & \text{if } i\equiv p,\\
    \{x_{i-1},\dots,x_{i+2}\} & \text{if } i\not\equiv p.
  \end{cases}
\]

\begin{remark}
For $m\ge5$ the three edges $\ell_0,\ell_3,\ell_{-3}$ are pairwise disjoint and any two of
them extend to a perfect matching, but the three together do not because $x_1$ cannot
be matched.
So $\{\ell_0,\ell_3,\ell_{-3}\}$ is a \emph{bad matching} in the sense of Chandrakar and Singh,
a minimal matching that does not extend to a perfect matching \cite[Def.~2.3]{PLTpm}, with three edges.
Hence $\Mp(\AP_m)$ is not a flag complex.
In contrast, the next lemma shows $K_0$ and $K_1$ are flag complexes.
\end{remark}

\begin{lemma}\label{lem:extend}
Let $p\in\{0,1\}$ be a parity.
A set $S$ of edges extends to a perfect matching of parity $p$ if and only if
\begin{enumerate}[label=(\roman*)]
\item every short edge $s_i\in S$ has $i\equiv p$, and
\item for any two edges $e,f\in S$, either $V_p(e)=V_p(f)$ or $V_p(e)\cap V_p(f)=\emptyset$.
\end{enumerate}
Consequently, $K_p$ is a flag complex.
\end{lemma}

\begin{proof}
If $S$ is contained in a perfect matching of parity $p$, then by the definition of parity the
sets $V_p(e)$, $e\in S$, belong to its partition in \cref{lem:tilings}, and (i) and (ii) follow.

Conversely, suppose that (i) and (ii) hold.
We prove the case $p=0$, and the case $p=1$ is similar.
By~(i) and the definition of $V_0$, every set $V_0(e)$, $e\in S$, consists of two or four
consecutive vertices beginning with an even index, and by~(ii) the distinct ones are pairwise
disjoint.
So the remaining vertices can be partitioned into pairs $\{x_{2k},x_{2k+1}\}$.
Together with the distinct sets $V_0(e)$, these pairs form a partition as in \cref{lem:tilings},
which determines a perfect matching of parity $0$ containing $S$.

Lastly, $K_p$ is a flag complex because conditions (i) and (ii) concern only single edges and pairs of edges of $S$.
\end{proof}

\begin{example}\label{ex:kp}
Let $m=4$ and $p=0$, and let $S=\{\ell_0,\ell_1,s_4\}$ (\cref{fig:graphs}).
Then $V_0(\ell_0)=V_0(\ell_1)=\{x_0,\dots,x_3\}$ and $V_0(s_4)=\{x_4,x_5\}$, and $S$ extends to
the perfect matching $\{\ell_0,\ell_1,s_4,s_6\}$ of parity $0$.
\end{example}

Now we compute the homotopy type of $K_p$. First we reduce the problem using elementary collapses. 
Recall that if $\tau$ and $\tau\cup\{x\}$ with $x\notin\tau$ are the only faces of a simplicial
complex $X$ that contain $\tau$, then $\tau$ is a \emph{free face} of $X$, and removing $\tau$ and
$\tau\cup\{x\}$ from $X$ is an \emph{elementary collapse}.
If a sequence of elementary collapses turns $X$ into $Y$, then $X$ \emph{collapses onto} $Y$,
and $X\simeq Y$ \cite[Sect.~6.4]{KozlovBook}.
In particular, if $a\ne b$ are elements of the ground set of $X$ such that every facet of $X$
containing $a$ also contains $b$, then $X$ collapses onto the subcomplex of faces not
containing $a$.
This is well known \cite[Remark~2.4]{BarmakMinian12}, and we include a short proof.
Every face $\sigma$ containing $a$ lies in a facet, which contains $b$, so $\sigma\cup\{b\}$ is a
face.
We delete the pairs $\sigma$, $\sigma\cup\{b\}$ with $a\in\sigma$ and $b\notin\sigma$ in the order
of decreasing $|\sigma|$.
When the pair of $\sigma$ is deleted, every face that properly contains $\sigma$, except
$\sigma\cup\{b\}$, lies in the pair of a larger set, and this pair is already deleted.
So $\sigma$ is a free face at that moment.
We use this in the proof of \cref{prop:Kp} to remove one of the two long edges in each pair
$\{\ell_{2k},\ell_{2k+1}\}$.

For an element $a$ of the ground set of $X$ and a face $\sigma$ of $X$, the \emph{star} of $a$
and the \emph{star cluster} of $\sigma$ \cite[Def.~3.1]{Barmak13} are the subcomplexes
\[
  \mathrm{st}_X(a)=\{\tau\in X:\tau\cup\{a\}\in X\},\qquad
  \mathrm{SC}_X(\sigma)=\bigcup_{a\in\sigma}\mathrm{st}_X(a).
\]

\begin{lemma}[{\cite[Lemma~3.2]{Barmak13}}]\label{lem:cluster}
The star cluster of a face of a flag complex is contractible.
\end{lemma}

\begin{proposition}\label{prop:Kp}
$K_p$ is contractible if $m$ is odd, and $K_p\simeq S^{m/2-1}\vee S^{m/2-1}$ if $m$ is even.
\end{proposition}

\begin{proof}
Since $K_0\cong K_1$, we may assume that $p=0$, and we write $K=K_0$ and $V=V_0$.
By \cref{lem:tilings}, for each $k$, the edges $\ell_{2k}$ and $\ell_{2k+1}$ lie in the same perfect
matchings of parity $0$.
So every facet of $K$ that contains $\ell_{2k+1}$ also contains $\ell_{2k}$, and $K$ collapses
onto its subcomplex of faces not containing $\ell_{2k+1}$.
Therefore $K$ collapses onto the subcomplex $K'$ of faces
that contain no long edges with odd indices.
Now the long edges $\ell_j$ in $K'$ have even indices, so $V(\ell_j)=\{x_j,\dots,x_{j+3}\}$.

Next, the short edges of $K'$ form the face $\Delta=\{s_0,s_2,\dots,s_{2m-2}\}$.
By \cref{lem:extend}, the faces of $K'$ are exactly the sets of edges that contain no $\ell_{2k+1}$ and
satisfy conditions (i) and (ii).
These conditions concern only single edges and pairs of edges, so $K'$ is a flag complex.
So the star cluster $A$ of $\Delta$ in $K'$ is contractible by \cref{lem:cluster}.
We claim that every face $\sigma$ of $K'$ satisfies
\begin{equation}\label{eq:cluster}
  \sigma\in A\iff\bigcup_{\ell_j\in\sigma}V(\ell_j)\ne\{x_0,x_1,\dots,x_{2m-1}\}.
\end{equation}
If $\sigma\in A$, then $\sigma\cup\{s_{2k}\}$ is a face of $K'$ for some $k$, and by
\cref{lem:extend} the sets $V(\ell_j)$, $\ell_j\in\sigma$, miss $V(s_{2k})=\{x_{2k},x_{2k+1}\}$.
Conversely, if these sets miss a vertex, then they miss a whole pair
$\{x_{2k},x_{2k+1}\}=V(s_{2k})$, because each of them is a union of such pairs.
Then $\sigma\cup\{s_{2k}\}$ is a face of $K'$ by \cref{lem:extend}, so $\sigma\in A$.

If $m$ is odd, then $K'=A$ by \eqref{eq:cluster}, since the union in \eqref{eq:cluster} is
disjoint by \cref{lem:extend}, so its size is a multiple of~$4$, while $2m$ is not.
Hence $K'$ is contractible.

On the other hand, if $m$ is even, then by \eqref{eq:cluster} and \cref{lem:extend} a face of
$K'$ is not in $A$ if and only if its sets $V(\ell_j)=\{x_j,\dots,x_{j+3}\}$ partition the vertex
set.
There are only two such partitions given by the faces
\[
  \Phi_0=\{\ell_0,\ell_4,\dots,\ell_{2m-4}\},\qquad
  \Phi_2=\{\ell_2,\ell_6,\dots,\ell_{2m-2}\}.
\]
Hence $K'=A\cup\{\Phi_0,\Phi_2\}$, so $K'$ is obtained from $A$ by attaching the two
$(m/2-1)$-simplices $\Phi_0$ and $\Phi_2$ along their boundaries, which lie in $A$.
Since $A$ is contractible,
\[
  K'\simeq K'/A\cong S^{m/2-1}\vee S^{m/2-1}.
  \qedhere
\]
\end{proof}

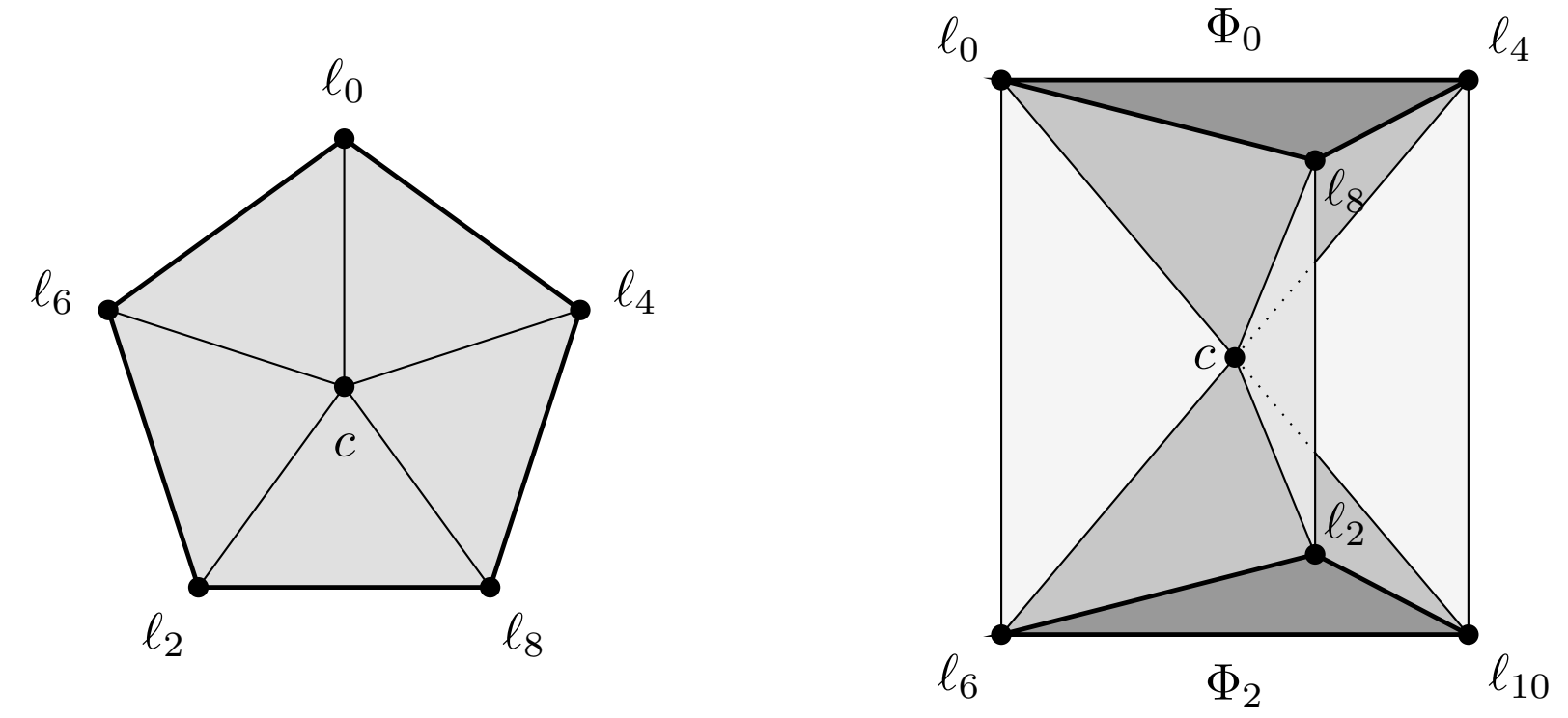
\begin{figure}[ht]
\centering
\begin{tikzpicture}[baseline=(current bounding box.center)]
  \coordinate (c) at (0,0);
  \foreach \a/\k in {0/0,4/1,8/2,2/3,6/4}{%
    \coordinate (p\a) at ({90-72*\k}:1.7);}
  \fill[black!12] (p0)--(p4)--(p8)--(p2)--(p6)--cycle;
  \draw[line width=0.9pt] (p0)--(p4)--(p8)--(p2)--(p6)--cycle;
  \foreach \a in {0,4,8,2,6}{\draw (c)--(p\a); \fill (p\a) circle (2pt);}
  \fill (c) circle (2pt);
  \node[font=\small] at (0,-0.4) {$c$};
  \foreach \a/\k in {0/0,4/1,8/2,2/3,6/4}{%
    \node[font=\small] at ({90-72*\k}:2.1) {$\ell_{\a}$};}
\end{tikzpicture}
\qquad\qquad
\begin{tikzpicture}[baseline=(current bounding box.center)]
  \coordinate (c) at (0,0);
  \coordinate (a0) at (-1.6,1.9); \coordinate (a4) at (1.6,1.9); \coordinate (a8) at (0.55,1.35);
  \coordinate (b6) at (-1.6,-1.9); \coordinate (b10) at (1.6,-1.9); \coordinate (b2) at (0.55,-1.35);
  % cones over the 1-simplices {l_0,l_6} and {l_4,l_10} (lightest)
  \fill[black!4] (c)--(a0)--(b6)--cycle;
  \fill[black!4] (c)--(a4)--(b10)--cycle;
  % cones over the boundaries of Phi_0 and Phi_2 (light grey)
  \fill[black!22] (c)--(a0)--(a4)--cycle;
  \fill[black!22] (c)--(b6)--(b10)--cycle;
  % the facets Phi_0 and Phi_2 (dark grey)
  \fill[black!40] (a0)--(a4)--(a8)--cycle;
  \fill[black!40] (b6)--(b10)--(b2)--cycle;
  \draw[line width=0.9pt] (a0)--(a4)--(a8)--cycle;
  \draw[line width=0.9pt] (b6)--(b10)--(b2)--cycle;
  \draw (a0)--(b6) (a4)--(b10);
  \foreach \q in {a0,a4,b6,b10}{\draw (c)--(\q);}
  % the cone over {l_8,l_2}, in front of the rest
  \fill[black!10] (c)--(a8)--(b2)--cycle;
  \begin{scope}
    \clip (c)--(a8)--(b2)--cycle;
    \draw[dotted] (c)--(a4) (c)--(b10);
  \end{scope}
  \draw (c)--(a8)--(b2)--cycle;
  \foreach \q in {a0,a4,a8,b6,b10,b2}{\fill (\q) circle (2pt);}
  \fill (c) circle (2pt);
  \node[font=\small,left] at (c) {$c$};
  \node[font=\small,above left] at (a0) {$\ell_0$};
  \node[font=\small,above right] at (a4) {$\ell_4$};
  \node[font=\small,below right,inner sep=1.5pt] at (a8) {$\ell_8$};
  \node[font=\small,below left] at (b6) {$\ell_6$};
  \node[font=\small,below right] at (b10) {$\ell_{10}$};
  \node[font=\small,above right,inner sep=1.5pt] at (b2) {$\ell_2$};
  \node[font=\small] at (0,2.25) {$\Phi_0$};
  \node[font=\small] at (0,-2.25) {$\Phi_2$};
\end{tikzpicture}
\caption{The complex $D\cup c*(A\cap D)\simeq K'$ for $\AP_5$ (left) and $\AP_6$
(right).}\label{fig:K1}
\end{figure}

Shrinking $A$ to a point gives a picture of $K'$.
Let $D$ be the subcomplex of faces of $K'$ that consist of long edges only.
Then $K'=A\cup D$, and $A\cap D$ consists of all faces of $D$ if $m$ is odd, and of all faces of
$D$ other than $\Phi_0$ and $\Phi_2$ if $m$ is even.
Let $c*(A\cap D)$ be the cone over $A\cap D$ with a new cone point $c$.
Since $A$ and $c*(A\cap D)$ are contractible, we have
\[
  K'\simeq K'/A=D/(A\cap D)\simeq D\cup c*(A\cap D).
\]
So the cone $c*(A\cap D)$ takes the place of $A$ (\cref{fig:K1}).

The next proposition gives the structure of the intersection of $K_0$ and $K_1$.
\begin{proposition}\label{prop:intersection}
The map $\ell_i\mapsto x_i$ is an isomorphism from $K_0\cap K_1$ onto $\Ind(C_{2m}^3)$.
\end{proposition}

\begin{proof}
By the definition of $C_{2m}^3$, the complex $\Ind(C_{2m}^3)$ is the flag complex whose faces
with at most two elements are the sets $\{x_i\}$ and the sets $\{x_i,x_k\}$ with
$k-i\notin\{\pm1,\pm2,\pm3\}$.
The complex $K_0\cap K_1$ is also a flag complex, since $K_0$ and $K_1$ are
(\cref{lem:extend}).
So it suffices to show that the faces of $K_0\cap K_1$ with at most two elements are the sets
$\{\ell_i\}$ and the sets $\{\ell_i,\ell_k\}$ with $k-i\notin\{\pm1,\pm2,\pm3\}$.

First, the faces of $K_0\cap K_1$ with one element are the sets $\{\ell_i\}$
since, by \cref{lem:extend}, $\{s_i\}$ lies only in $K_p$ with $p\equiv i$.
Consequently, every face of $K_0\cap K_1$ with two elements has the form $\{\ell_i,\ell_k\}$.
By rotational symmetry of the graph, we may assume $i=0$ and $1\le k\le m$, and we show that
$\{\ell_0,\ell_k\}$ is a face of $K_0\cap K_1$ if and only if $k\ge4$.
By \cref{lem:extend}, $\{\ell_0,\ell_k\}$ is a face of $K_0\cap K_1$ if and only if for both
$p=0,1$ the sets $V_p(\ell_0)$ and $V_p(\ell_k)$ are equal or disjoint.
\begin{itemize}
\item If $k=1$, then $V_1(\ell_0)=\{x_{-1},\dots,x_2\}$ and $V_1(\ell_1)=\{x_1,\dots,x_4\}$
  share $x_1$ and $x_2$.
\item If $k=2$, then $V_0(\ell_0)=\{x_0,\dots,x_3\}$ and $V_0(\ell_2)=\{x_2,\dots,x_5\}$
  share $x_2$ and $x_3$.
\item If $k=3$, then $V_0(\ell_0)=\{x_0,\dots,x_3\}$ and $V_0(\ell_3)=\{x_2,\dots,x_5\}$
  share $x_2$ and $x_3$.
\item If $k\ge4$ is even, then $V_0(\ell_0)$ and $V_0(\ell_k)$ begin at $x_0$ and $x_k$, and
  $V_1(\ell_0)$ and $V_1(\ell_k)$ begin at $x_{-1}$ and $x_{k-1}$.
  Both differences are $k$, and $4\le k\le m\le 2m-4$, so the condition holds.
\item If $k\ge5$ is odd, then $V_0(\ell_0)$ and $V_0(\ell_k)$ begin at $x_0$ and $x_{k-1}$, and
  $V_1(\ell_0)$ and $V_1(\ell_k)$ begin at $x_{-1}$ and $x_k$.
  The differences are $k-1$ and $k+1$, and $4\le k-1<k+1\le m+1\le 2m-4$ since $m\ge k\ge5$,
  so the condition holds.
\end{itemize}
\end{proof}

The following well-known lemma is used to glue $K_0$ and $K_1$. We include a proof for completeness.

\begin{lemma}\label{lem:glue}
Let $X$ be a simplicial complex that is the union of two connected subcomplexes $Y$ and $Z$
whose intersection $W=Y\cap Z$ is nonempty.
If the inclusions $W\to Y$ and $W\to Z$ are null-homotopic, then $X\simeq Y\vee Z\vee\Sigma W$.
In particular, if $Z$ is contractible and the inclusion $W\to Y$ is null-homotopic, then
$X\simeq Y\vee\Sigma W$, and if $Y$ and $Z$ are contractible, then $X\simeq\Sigma W$.
\end{lemma}

\begin{proof}
Let $M$ be obtained from $Y\sqcup(W\times[0,1])\sqcup Z$ by gluing $W\times\{0\}$ to
$W\subseteq Y$ and $W\times\{1\}$ to $W\subseteq Z$.
The projection $M\to X$, $(w,t)\mapsto w$, is a homotopy equivalence, since $W$ is a
subcomplex of $Y$ and of $Z$.
On the other hand, $M$ is obtained from $Y\sqcup Z$ by attaching $W\times[0,1]$ along
$W\times\{0,1\}$.
Replacing the two null-homotopic attaching maps by constant maps turns the cylinder into $\Sigma W$, attached to $Y$ and $Z$ at its two cone points, and this
space is homotopy equivalent to $Y\vee Z\vee\Sigma W$.
\end{proof}

\begin{theorem}\label{thm:main}
Let $m\ge3$. Then
\[
  \Mp(\AP_m)\simeq
  \begin{cases}
    \Sigma\Ind(C_{2m}^3) & \text{if $m$ is odd},\\[2pt]
    \bigvee_4 S^{m/2-1}\ \vee\ \Sigma\Ind(C_{2m}^3) & \text{if $m$ is even}.
  \end{cases}
\]
\end{theorem}

\begin{proof}
Let $I=K_0\cap K_1\cong\Ind(C_{2m}^3)$ (\cref{prop:intersection}).
Since $\Mp(\AP_m)=K_0\cup K_1$, by \cref{lem:glue,prop:Kp} it suffices to show that the
inclusion $I\to K_p$ is null-homotopic, and for odd $m$ this is clear since $K_p$ is
contractible.
Let $m$ be even, and assume $p=0$ without loss of generality. 
The case $p=1$ follows from index shifting $i\rightarrow i+1$.

The vertices of an independent set of $C_{2m}^3$ are pairwise at distance at least $4$ along
the $2m$-cycle, so it has at most $m/2$ vertices.
The only independent sets with $m/2$ vertices are the sets $\{x_c,x_{c+4},\dots,x_{c+2m-4}\}$,
$c\in\{0,1,2,3\}$, so the only faces of $I$ with $m/2$ elements are the four sets
$\sigma_c=\{\ell_c,\ell_{c+4},\dots,\ell_{c+2m-4}\}$.
The face $\sigma_c\setminus\{\ell_c\}$ contains $\ell_{c+4}$, which lies in no $\sigma_{c'}$ other
than $\sigma_c$, so it is a free face.
Hence $I$ collapses onto the subcomplex $I'$ obtained by deleting $\sigma_c$ and
$\sigma_c\setminus\{\ell_c\}$ for all $c$.

Let $B$ be the star cluster of $\Delta=\{s_0,s_2,\dots,s_{2m-2}\}$ in $K_0$.
It is contractible by \cref{lem:cluster}, since $K_0$ is a flag complex by \cref{lem:extend}.
As in the proof of \cref{prop:Kp}, every face $\sigma$ of $K_0$ satisfies
\begin{equation}\label{eq:cluster0}
  \sigma\in B\iff\bigcup_{\ell_i\in\sigma}V_0(\ell_i)\ne\{x_0,\dots,x_{2m-1}\},
\end{equation}
since every $V_0(\ell_i)$ begins at an even index, so the proof of \eqref{eq:cluster} applies.
We show that every face $\sigma\in I\setminus B$ is one of the sets $\sigma_c$, so that $I'\subseteq B$.
By \eqref{eq:cluster0}, the sets $V_0(\ell_i)$, $\ell_i\in\sigma$, cover the whole vertex set.
Each of them has four vertices, so $\sigma$ has at least $m/2$ elements.
Since the faces of $I$ have at most $m/2$ elements, $\sigma$ is one of the sets $\sigma_c$.
Therefore the inclusion $I\to K_0$ is null-homotopic, since $I\simeq I'\subseteq B$ and $B$ is
contractible.
\end{proof}

\section{The explicit homotopy type}\label{sec:explicit}

By \cref{thm:main}, $\Mp(\AP_m)$ is determined by $\Ind(C_{2m}^3)$.
Adamaszek gave a recursion for the independence complexes of powers of cycles
\cite[Thm.~1.1 and Cor.~6.5]{Adamaszek}. For the cube and $n\ge19$ it reads
\begin{equation}\label{eq:cycle-rec}
  \begin{aligned}
    \Ind(C_n^3)\simeq{}&\Sigma^2\Ind(C_{n-12}^3)\\
    &\vee\Sigma^3\Bigl(\bigvee_5\Ind(P_{n-18}^3)\vee\bigvee_{11}\Ind(P_{n-19}^3)\\
    &\qquad\vee\bigvee_{13}\Ind(P_{n-20}^3)\vee\bigvee_7\Ind(P_{n-21}^3)\Bigr),
  \end{aligned}
\end{equation}
where the powers of paths satisfy, for $n\ge4$,
\begin{equation}\label{eq:path-rec}
  \Ind(P_n^3)\simeq\Sigma\Ind(P_{n-5}^3)\vee\Sigma\Ind(P_{n-6}^3)\vee\Sigma\Ind(P_{n-7}^3)
\end{equation}
\cite[(3)]{Adamaszek}, \cite{Engstrom09}, with $\Ind(P_n^3)=S^{-1}$ for $-3\le n\le0$ (the
empty graph), and $\Ind(P_n^3)$ a set of $n$ points for $1\le n\le3$, since any two of the $n$
vertices of $P_n^3$ are adjacent.
The recursion \eqref{eq:cycle-rec} says nothing about $n<19$ \cite[Remark~6.6]{Adamaszek}, and
for $2m\ge19$ it reduces $\Ind(C_{2m}^3)$ to $\Ind(C_{2m-12}^3)$.
So the initial values $\Ind(C_n^3)$, $n\le18$, are needed (for $\Mp(\AP_m)$ only the even
ones).

We compute the initial values using discrete Morse theory in the form of \cite[Ch.~4]{Jonsson}.
Let $G$ be a graph whose vertices are ordered $v_1,\dots,v_m$.
We include the empty set, of dimension $-1$, among the faces of $\Ind(G)$.
Put $S_0=\Ind(G)$.
For $j=1,\dots,m$, the \emph{step of $v_j$} pairs $\sigma$ with $\sigma\cup\{v_j\}$
whenever $v_j\notin\sigma$ and both $\sigma$ and $\sigma\cup\{v_j\}$ lie in $S_{j-1}$, and $S_j$
is $S_{j-1}$ without the faces paired in this step.
After the last step, the faces in $S_m$ are \emph{critical}.
By \cite[Lemma~4.1]{Jonsson} and induction on $m$, all these pairs together form an acyclic
matching in the sense of discrete Morse theory, in which the empty set is paired with $\{v_1\}$.
Hence $\Ind(G)$ is homotopy equivalent to a CW complex with one $0$-cell and one $d$-cell for
each critical face of dimension $d$ \cite[Thm.~6.3]{FormanGuide}, \cite[Thm.~11.13(b)]{KozlovBook}.
If all critical faces have the same dimension $d$, then $\Ind(G)$ is a wedge of $d$-spheres,
one for each critical face.

\begin{proposition}\label{prop:initial}
For $6\le n\le18$, $\Ind(C_n^3)$ is the wedge of spheres in \cref{tab:initial}.
\end{proposition}

\begin{proof}
We start from the homology, which can be computed with a computer algebra system.
A computation with SageMath \cite{Sage} shows that for $6\le n\le18$,
$\tilde H_*(\Ind(C_n^3);\Z)$ is free, of the ranks given by the wedges in \cref{tab:initial}. For instance, the only nonzero reduced homology groups of $\Ind(C_{15}^3)$
are $\tilde H_1\cong\mathbb Z$ and $\tilde H_2\cong\mathbb Z^5$.
Using a computer program, we loop over vertex orders until the number of critical faces in each
dimension equals the rank of the reduced homology.

For $n\ne14,15$, the order $x_0,x_1,\dots,x_{n-1}$ works, and all its critical faces have the
same dimension, so $\Ind(C_n^3)$ is the wedge of spheres in \cref{tab:initial}.

For $n=14,15$, the order $x_0,x_1,\dots,x_{n-1}$ leaves more critical faces than the ranks of
the homology, but the orders
\[
  \begin{aligned}
    &x_1,x_5,x_{11},x_3,x_9,x_7,x_2,x_{13},x_0,x_{10},x_8,x_4,x_{12},x_6 && (n=14),\\
    &x_0,x_{11},x_7,x_{14},x_{10},x_8,x_{13},x_3,x_{12},x_9,x_5,x_4,x_6,x_2,x_1 && (n=15)
  \end{aligned}
\]
leave numbers of critical faces matching the ranks.
Hences $\Ind(C_{14}^3)\simeq\bigvee_8S^1$, and $\Ind(C_{15}^3)$ is homotopy equivalent to a CW complex $X$ with one $0$-cell, one
$1$-cell and five $2$-cells.
The $1$-skeleton of $X$ is a circle, and the $k$th $2$-cell is attached by a map
$S^1\to S^1$ of some degree $d_k$, so $H_1(X;\mathbb Z)\cong\mathbb Z/(d_1,\dots,d_5)$.
Since $H_1(\Ind(C_{15}^3);\mathbb Z)\cong\mathbb Z$, we get $d_1=\dots=d_5=0$, and the
attaching maps are null-homotopic.
Hence $\Ind(C_{15}^3)\simeq S^1\vee\bigvee_5S^2$.
\end{proof}

\begin{corollary}\label{cor:wedge}
For every $n\ge6$, $\Ind(C_n^3)$ is homotopy equivalent to a wedge of spheres, and so is
$\Mp(\AP_m)$ for every $m\ge3$.
The numbers of spheres are determined by \eqref{eq:cycle-rec}, \eqref{eq:path-rec},
\cref{tab:initial} and \cref{thm:main}.
\end{corollary}

\begin{table}[ht]
\centering
\renewcommand{\arraystretch}{1.2}
\begin{tabular}{@{}clll@{}}
  \toprule
  $n$ & $\Ind(P_n^3)$ & $\Ind(C_n^3)$ & $\Mp(\AP_{n/2})$\\
  \midrule
  $1$ & point &  & \\
  $2$ & $S^0$ &  & \\
  $3$ & $\bigvee_2S^0$ &  & \\
  $4$ & $\bigvee_3S^0$ &  & \\
  $5$ & $\bigvee_3S^0$ &  & \\
  $6$ & $\bigvee_2S^0$ & $\bigvee_5S^0$ & $\bigvee_5S^1$\\
  $7$ & $S^0\vee S^1$ & $\bigvee_6S^0$ & \\
  $8$ & $\bigvee_3S^1$ & $\bigvee_3S^0$ & $\bigvee_7S^1$\\
  $9$ & $\bigvee_6S^1$ & $S^1$ & \\
  $10$ & $\bigvee_8S^1$ & $\bigvee_6S^1$ & $\bigvee_6S^2$\\
  $11$ & $\bigvee_8S^1$ & $\bigvee_{12}S^1$ & \\
  $12$ & $\bigvee_6S^1\vee S^2$ & $\bigvee_{15}S^1$ & $\bigvee_{19}S^2$\\
  $13$ & $\bigvee_3S^1\vee\bigvee_4S^2$ & $\bigvee_{14}S^1$ & \\
  $14$ &  & $\bigvee_8S^1$ & $\bigvee_8S^2$\\
  $15$ &  & $S^1\vee\bigvee_5S^2$ & \\
  $16$ &  & $\bigvee_{19}S^2$ & $\bigvee_{23}S^3$\\
  $17$ &  & $\bigvee_{33}S^2$ & \\
  $18$ &  & $\bigvee_{41}S^2$ & $\bigvee_{41}S^3$\\
  \bottomrule
\end{tabular}
\caption{Initial values of $\Ind(P_n^3)$, $\Ind(C_n^3)$, and
$\Mp(\AP_m)$, $n=2m$.}
\label{tab:initial}
\end{table}

\section{The independence complex of \texorpdfstring{$\GP(n,2)$}{GP(n,2)}}\label{sec:gp}

In this section we study the homotopy type of $\Ind(\GP(n,2))$ for all $n\ge5$.
We first recall two tools for independence complexes.
Let $z$ be a vertex of a graph $H$.
There is a cofibre sequence (see \cite[Sect.~1]{Adamaszek})
\[
  \Ind(H\setminus N[z])\xhookrightarrow{\;i\;}\Ind(H\setminus z)\hookrightarrow\Ind(H),
\]
where $i$ is the inclusion.
Here $\Ind(H\setminus z)$ and $\Ind(H\setminus N[z])$ are the deletion and the link of the
$0$-simplex $\{z\}$ in $\Ind(H)$.
The following lemma is an immediate consequence of the cofibre sequence.

\begin{lemma}[{see \cite[Prop.~3.1]{Adamaszek}}]\label{lem:cofibre}
Let $z$ be a vertex of a graph $H$.
\begin{enumerate}[label=(\alph*)]
\item If $\Ind(H\setminus N[z])$ is contractible, then $\Ind(H)\simeq\Ind(H\setminus z)$.
\item If $\Ind(H\setminus z)$ is contractible, then $\Ind(H)\simeq\Sigma\Ind(H\setminus N[z])$.
\item If $i$ is null-homotopic, then
  $\Ind(H)\simeq\Ind(H\setminus z)\vee\Sigma\Ind(H\setminus N[z])$.
  This is the case if, for some $k\ge0$, $\Ind(H\setminus N[z])$ is homotopy equivalent to a
  wedge of $k$-spheres and $\Ind(H\setminus z)$ to a wedge of $(k+1)$-spheres.
\end{enumerate}
\end{lemma}

\begin{lemma}[Fold lemma {\cite[Lemma~3.2]{Engstrom09}}]\label{lem:fold}
Let $z$ and $w$ be distinct vertices of a graph $H$ with $N(z)\subseteq N(w)$.
Then $\Ind(H)$ collapses onto $\Ind(H\setminus w)$.
\end{lemma}

The following consequence of the fold lemma is used in the remainder of the paper repeatedly.
It is the case of a single neighbour in \cite[Thm.~3.7]{Engstrom09} (see also
\cite[Cor.~5.4]{Adamaszek}), and for the path $P_n$ with $n$ vertices it is Kozlov's recursion
$\Ind(P_n)\simeq\Sigma\Ind(P_{n-3})$ \cite[Prop.~4.6]{Kozlov99}.

\begin{lemma}[{\cite[Thm.~3.7]{Engstrom09}}]\label{lem:split}
Let $y$ be a vertex of a graph $H$, and let $x$ be a neighbour of $y$ whose only neighbour
is $y$.
Then $\Ind(H)\simeq\Sigma\Ind(H\setminus N[y])$.
\end{lemma}

\begin{proof}
For every neighbour $w\ne x$ of $y$ we have $N(x)=\{y\}\subseteq N(w)$, so the fold lemma
deletes $w$.
Then $xy$ is a connected component, and
$\Ind((H\setminus N[y])\sqcup K_2)=\Sigma\Ind(H\setminus N[y])$.
\end{proof}

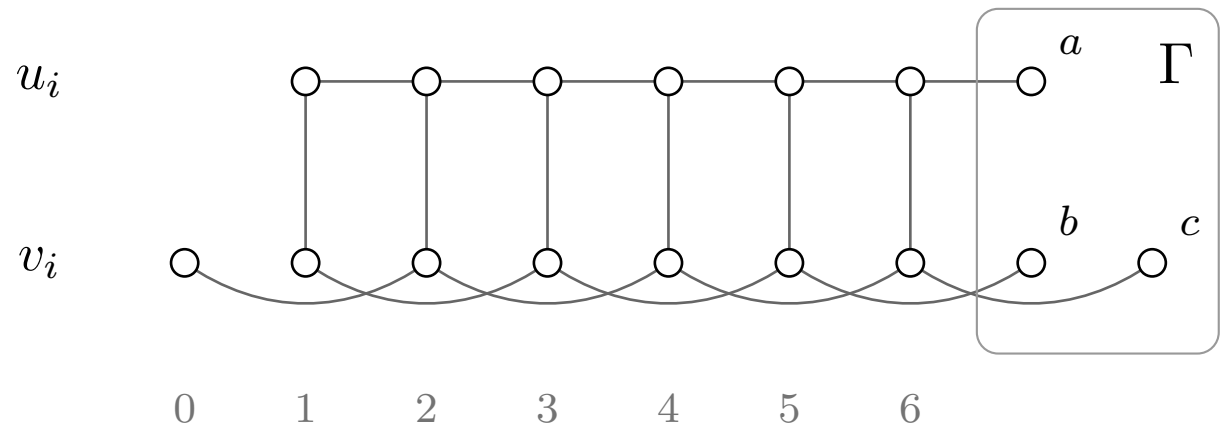
\begin{figure}[ht]
\centering
\begin{tikzpicture}[x=8mm,y=12mm]
  \draw[gpedge] (6,1) -- (7,1);
  \foreach \i in {5,6} { \draw[gpinner] (\i,0) to (\i+2,0); }
  \draw[rounded corners,black!40] (6.55,-0.5) rectangle (8.55,1.4);
  \node at (8.2,1.1) {$\Gamma$};
  \gpgraph{0,...,6}{1,...,6}{0,...,6}
  \foreach \x/\y/\t in {7/1/a,7/0/b,8/0/c} {
    \node[gpvertex] at (\x,\y) {};
    \node[font=\scriptsize,anchor=south west] at (\x+0.06,\y+0.04) {$\t$}; }
  \gpcols[-0.8]{0,...,6}
  \node[font=\small] at (-1.2,1) {$u_i$};
  \node[font=\small] at (-1.2,0) {$v_i$};
\end{tikzpicture}
\caption{The graph $R_6(\Gamma)$. Column $0$ has lost $u_0$, and the three edges $u_6a$,
$v_5b$, $v_6c$ cross from column $6$ into $\Gamma$.}\label{fig:strip}
\end{figure}

Cutting $\GP(n,2)$ open near one column leaves graphs of the following form.
We draw $\GP(n,2)$ by \emph{columns} $\{u_i,v_i\}$, side by side (\cref{fig:strip}).
For $m\ge0$ let $S_m$ be the graph with vertices $v_0$ and $u_i,v_i$ for $1\le i\le m$, and
with the edges $u_iu_{i+1}$, $u_iv_i$, $v_iv_{i+2}$ among them.
Thus $S_m$ consists of the columns $0,\dots,m$, where column $0$ has lost $u_0$.
Let $\Gamma$ be a graph and let $a,b,c$ be vertices of $\Gamma$, any of which may also be
absent.
Let $R_m(\Gamma)$ be the graph obtained from the disjoint union of $S_m$ and $\Gamma$ by adding
the three edges
\[
  u_ma,\qquad v_{m-1}b,\qquad v_mc
\]
that cross from column $m$ to the right.
An edge is omitted if one of its endpoints is absent. For $m=0$ only $v_0c$ is added.

\begin{proposition}\label{prop:block}
Let $\Gamma$ be as above and $m\ge5$. Then
\[
  \Ind(R_m(\Gamma))\simeq\Sigma^3\Ind(R_{m-5}(\Gamma)).
\]
\end{proposition}

\begin{proof}
Since $m\ge5$, the vertices $v_0,\dots,v_3$ and $u_1,\dots,u_4$ have no neighbours in
$\Gamma$.
We apply \cref{lem:split} three times (see \cref{fig:block}).
\begin{enumerate}[label=(\arabic*)]
\item The vertex $v_0$ has the single neighbour $v_2$. Applying \cref{lem:split} to $v_2$
  deletes $v_0,v_2,u_2,v_4$.
\item Now $u_1$ has the single neighbour $v_1$. Applying \cref{lem:split} to $v_1$ deletes
  $u_1,v_1,v_3$.
\item Now $u_3$ has the single neighbour $u_4$. Applying \cref{lem:split} to $u_4$ deletes
  $u_3,u_4,u_5$.
\end{enumerate}
This deletes the columns $0,\dots,4$ and the vertex $u_5$.
What is left is $R_{m-5}(\Gamma)$.
\end{proof}

\begin{figure}[ht]
\centering
\begin{tikzpicture}[x=5.5mm,y=11mm,baseline=(current bounding box.north)]
  \gpgraph[right]{0,...,8}{1,...,8}{0,...,8}
  \node[gpgrey] at (0,0) {};
  \node[gpgrey] at (2,1) {};
  \node[gpgrey] at (4,0) {};
  \node[gpblack] at (2,0) {};
  \node[font=\scriptsize,anchor=south west] at (0.06,0.04) {$x$};
  \node[font=\scriptsize,anchor=south west] at (2.06,0.04) {$y$};
  \gpcols[-0.8]{0,...,8}
  \node[font=\small,anchor=west] at (-0.4,1.55) {(1) $R_m(\Gamma)$};
\end{tikzpicture}
\hfill
\begin{tikzpicture}[x=5.5mm,y=11mm,baseline=(current bounding box.north)]
  \gpgraph[right]{0,...,8}{1,3,4,...,8}{1,3,5,6,7,8}
  \node[gpgrey] at (1,1) {};
  \node[gpgrey] at (3,0) {};
  \node[gpblack] at (1,0) {};
  \node[font=\scriptsize,anchor=south] at (1,1.08) {$x$};
  \node[font=\scriptsize,anchor=south west] at (1.06,0.04) {$y$};
  \gpcols[-0.8]{0,...,8}
  \node[font=\small,anchor=west] at (-0.4,1.55) {(2)};
\end{tikzpicture}

\bigskip
\begin{tikzpicture}[x=5.5mm,y=11mm,baseline=(current bounding box.north)]
  \gpgraph[right]{0,...,8}{3,...,8}{5,...,8}
  \node[gpgrey] at (3,1) {};
  \node[gpgrey] at (5,1) {};
  \node[gpblack] at (4,1) {};
  \node[font=\scriptsize,anchor=south] at (3,1.08) {$x$};
  \node[font=\scriptsize,anchor=south] at (4,1.08) {$y$};
  \gpcols[-0.8]{0,...,8}
  \node[font=\small,anchor=west] at (-0.4,1.55) {(3)};
\end{tikzpicture}
\hfill
\begin{tikzpicture}[x=5.5mm,y=11mm,baseline=(current bounding box.north)]
  \gpgraph[right]{0,...,8}{6,7,8}{5,...,8}
  \gpcols[-0.8]{0,...,8}
  \node[font=\small,anchor=west] at (-0.4,1.55) {(4) $R_{m-5}(\Gamma)$};
\end{tikzpicture}
\caption{The proof of \cref{prop:block}.}\label{fig:block}
\end{figure}

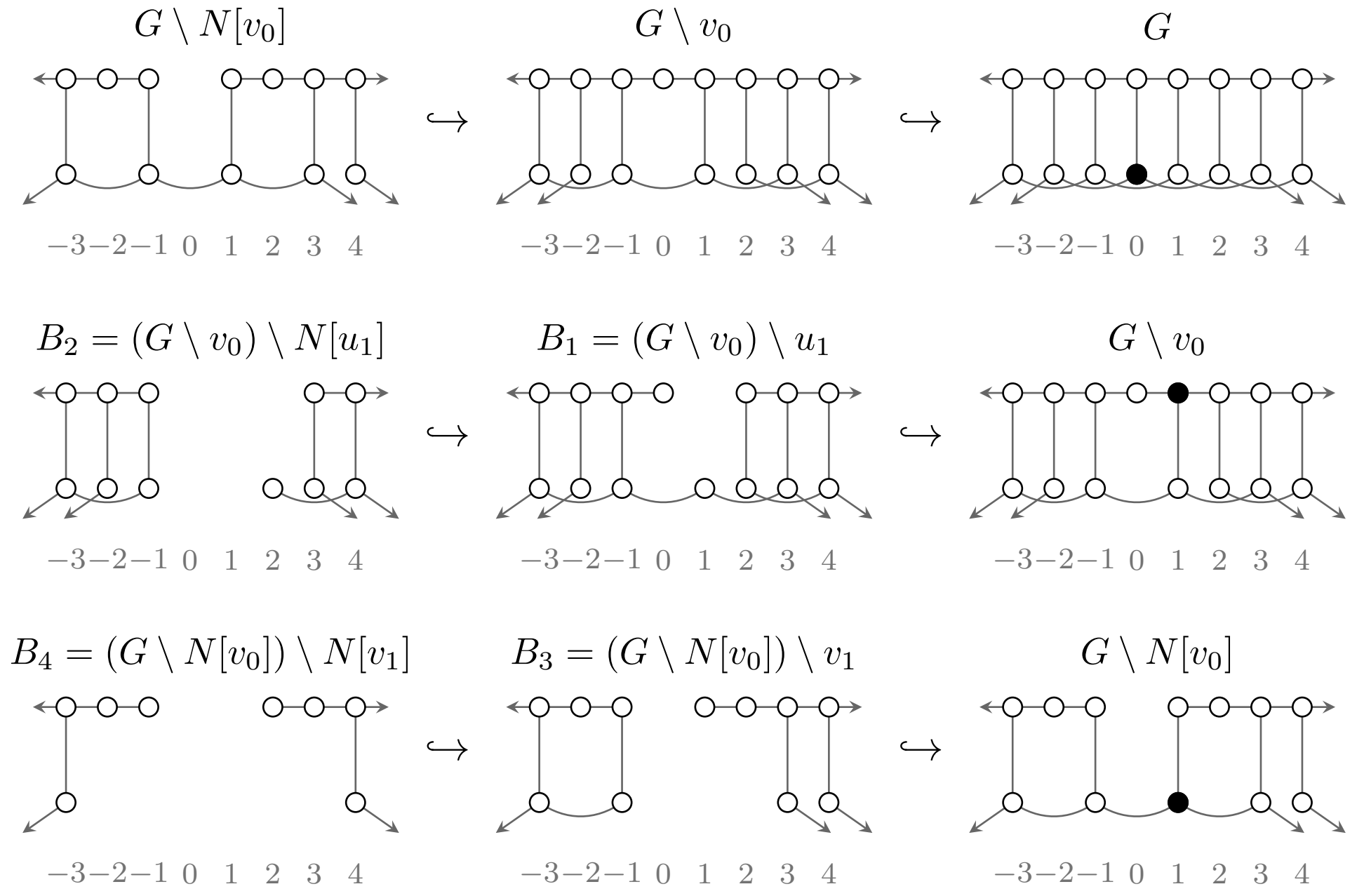
\begin{figure}[ht]
\centering
\begin{tikzpicture}[x=3.9mm,y=9mm,baseline={(0,0.4)}]
  \gpgraph[left,right]{-3,...,4}{-3,-2,-1,1,2,3,4}{-3,-1,1,3,4}
  \gpcols{-3,...,4}
  \node[font=\small] at (0.5,1.55) {$G\setminus N[v_0]$};
\end{tikzpicture}
\;$\hookrightarrow$\;
\begin{tikzpicture}[x=3.9mm,y=9mm,baseline={(0,0.4)}]
  \gpgraph[left,right]{-3,...,4}{-3,-2,-1,0,1,2,3,4}{-3,-2,-1,1,2,3,4}
  \gpcols{-3,...,4}
  \node[font=\small] at (0.5,1.55) {$G\setminus v_0$};
\end{tikzpicture}
\;$\hookrightarrow$\;
\begin{tikzpicture}[x=3.9mm,y=9mm,baseline={(0,0.4)}]
  \gpgraph[left,right]{-3,...,4}{-3,-2,-1,0,1,2,3,4}{-3,-2,-1,0,1,2,3,4}
  \node[gpblack] at (0,0) {};
  \gpcols{-3,...,4}
  \node[font=\small] at (0.5,1.55) {$G$};
\end{tikzpicture}

\bigskip
\begin{tikzpicture}[x=3.9mm,y=9mm,baseline={(0,0.4)}]
  \gpgraph[left,right]{-3,...,4}{-3,-2,-1,3,4}{-3,-2,-1,2,3,4}
  \gpcols{-3,...,4}
  \node[font=\small] at (0.5,1.55) {\makebox[0pt]{$B_2=(G\setminus v_0)\setminus N[u_1]$}};
\end{tikzpicture}
\;$\hookrightarrow$\;
\begin{tikzpicture}[x=3.9mm,y=9mm,baseline={(0,0.4)}]
  \gpgraph[left,right]{-3,...,4}{-3,-2,-1,0,2,3,4}{-3,-2,-1,1,2,3,4}
  \gpcols{-3,...,4}
  \node[font=\small] at (0.5,1.55) {\makebox[0pt]{$B_1=(G\setminus v_0)\setminus u_1$}};
\end{tikzpicture}
\;$\hookrightarrow$\;
\begin{tikzpicture}[x=3.9mm,y=9mm,baseline={(0,0.4)}]
  \gpgraph[left,right]{-3,...,4}{-3,-2,-1,0,1,2,3,4}{-3,-2,-1,1,2,3,4}
  \node[gpblack] at (1,1) {};
  \gpcols{-3,...,4}
  \node[font=\small] at (0.5,1.55) {$G\setminus v_0$};
\end{tikzpicture}

\bigskip
\begin{tikzpicture}[x=3.9mm,y=9mm,baseline={(0,0.4)}]
  \gpgraph[left,right]{-3,...,4}{-3,-2,-1,2,3,4}{-3,4}
  \gpcols{-3,...,4}
  \node[font=\small] at (0.5,1.55) {\makebox[0pt]{$B_4=(G\setminus N[v_0])\setminus N[v_1]$}};
\end{tikzpicture}
\;$\hookrightarrow$\;
\begin{tikzpicture}[x=3.9mm,y=9mm,baseline={(0,0.4)}]
  \gpgraph[left,right]{-3,...,4}{-3,-2,-1,1,2,3,4}{-3,-1,3,4}
  \gpcols{-3,...,4}
  \node[font=\small] at (0.5,1.55) {\makebox[0pt]{$B_3=(G\setminus N[v_0])\setminus v_1$}};
\end{tikzpicture}
\;$\hookrightarrow$\;
\begin{tikzpicture}[x=3.9mm,y=9mm,baseline={(0,0.4)}]
  \gpgraph[left,right]{-3,...,4}{-3,-2,-1,1,2,3,4}{-3,-1,1,3,4}
  \node[gpblack] at (1,0) {};
  \gpcols{-3,...,4}
  \node[font=\small] at (0.5,1.55) {$G\setminus N[v_0]$};
\end{tikzpicture}
\caption{The graphs of the three cofibre sequences for $\GP(n,2)$.
The graphs continue on both sides.}\label{fig:cofibre}
\end{figure}

From now on, $G=\GP(n,2)$ with $n\ge5$, and indices are in $\Z_n$.
We apply the cofibre sequence at $v_0$ in $G$, then at $u_1$ in $G\setminus v_0$ and at $v_1$ in
$G\setminus N[v_0]$.
This gives four graphs:
\[
  \begin{aligned}
    B_1&=G\setminus\{v_0,u_1\}, & B_2&=(G\setminus v_0)\setminus N[u_1],\\
    B_3&=(G\setminus N[v_0])\setminus v_1, & B_4&=(G\setminus N[v_0])\setminus N[v_1].
  \end{aligned}
\]
\Cref{fig:cofibre} shows the graphs of the three cofibre sequences.
After one application of \cref{lem:split} to $B_1$, $B_3$ and $B_4$, each $B_j$ is a graph
$R_m(\Gamma_j)$ as in \cref{fig:strip}.
Its strip starts to the right of the removed vertices, runs around the graph, 
and ends on the left with a small graph $\Gamma_j$, which has a different shape for each $j$.

\begin{figure}[ht]
\centering
\setlength{\tabcolsep}{3pt}
\begin{tabular}{@{}cccccc@{}}
 & $m=0$ & $m=1$ & $m=2$ & $m=3$ & $m=4$\\[4pt]
$\Gamma_1$
&\begin{tikzpicture}[x=4.4mm,y=8mm,baseline=(current bounding box.center)]
  \path (0,1.35);
  \draw[rounded corners,black!40] (0.55,-0.4) rectangle (1.45,1.3);
  \gpgraph{0,...,1}{}{0,1}
  \gpcols{0,...,1}
  \node[font=\small] at (0.5,-1.4) {$*$};
\end{tikzpicture}
&\begin{tikzpicture}[x=4.4mm,y=8mm,baseline=(current bounding box.center)]
  \path (0,1.35);
  \draw[rounded corners,black!40] (1.55,-0.4) rectangle (2.45,1.3);
  \gpgraph{0,...,2}{1}{0,1,2}
  \gpcols{0,...,2}
  \node[font=\small] at (1.0,-1.4) {$S^{1}$};
\end{tikzpicture}
&\begin{tikzpicture}[x=4.4mm,y=8mm,baseline=(current bounding box.center)]
  \path (0,1.35);
  \draw[rounded corners,black!40] (2.55,-0.4) rectangle (3.45,1.3);
  \gpgraph{0,...,3}{1,2}{0,1,2,3}
  \gpcols{0,...,3}
  \node[font=\small] at (1.5,-1.4) {$S^{1}$};
\end{tikzpicture}
&\begin{tikzpicture}[x=4.4mm,y=8mm,baseline=(current bounding box.center)]
  \path (0,1.35);
  \draw[rounded corners,black!40] (3.55,-0.4) rectangle (4.45,1.3);
  \gpgraph{0,...,4}{1,2,3}{0,1,2,3,4}
  \gpcols{0,...,4}
  \node[font=\small] at (2.0,-1.4) {$*$};
\end{tikzpicture}
&\begin{tikzpicture}[x=4.4mm,y=8mm,baseline=(current bounding box.center)]
  \path (0,1.35);
  \draw[rounded corners,black!40] (4.55,-0.4) rectangle (5.45,1.3);
  \gpgraph{0,...,5}{1,2,3,4}{0,1,2,3,4,5}
  \gpcols{0,...,5}
  \node[font=\small] at (2.5,-1.4) {$*$};
\end{tikzpicture}\\[8pt]
$\Gamma_2$
&\begin{tikzpicture}[x=4.4mm,y=8mm,baseline=(current bounding box.center)]
  \path (0,1.35);
  \gpgraph{0,...,0}{}{0}
  \gpcols{0,...,0}
  \node[font=\small] at (0.0,-1.4) {$*$};
\end{tikzpicture}
&\begin{tikzpicture}[x=4.4mm,y=8mm,baseline=(current bounding box.center)]
  \path (0,1.35);
  \gpgraph{0,...,1}{1}{0,1}
  \gpcols{0,...,1}
  \node[font=\small] at (0.5,-1.4) {$*$};
\end{tikzpicture}
&\begin{tikzpicture}[x=4.4mm,y=8mm,baseline=(current bounding box.center)]
  \path (0,1.35);
  \gpgraph{0,...,2}{1,2}{0,1,2}
  \gpcols{0,...,2}
  \node[font=\small] at (1.0,-1.4) {$S^{1}$};
\end{tikzpicture}
&\begin{tikzpicture}[x=4.4mm,y=8mm,baseline=(current bounding box.center)]
  \path (0,1.35);
  \gpgraph{0,...,3}{1,2,3}{0,1,2,3}
  \gpcols{0,...,3}
  \node[font=\small] at (1.5,-1.4) {$*$};
\end{tikzpicture}
&\begin{tikzpicture}[x=4.4mm,y=8mm,baseline=(current bounding box.center)]
  \path (0,1.35);
  \gpgraph{0,...,4}{1,2,3,4}{0,1,2,3,4}
  \gpcols{0,...,4}
  \node[font=\small] at (2.0,-1.4) {$S^{2}$};
\end{tikzpicture}\\[8pt]
$\Gamma_3$
&\begin{tikzpicture}[x=4.4mm,y=8mm,baseline=(current bounding box.center)]
  \path (0,1.35);
  \draw[rounded corners,black!40] (0.55,-0.4) rectangle (2.45,1.3);
  \gpgraph{0,...,2}{1,2}{0,2}
  \gpcols{0,...,2}
  \node[font=\small] at (1.0,-1.4) {$*$};
\end{tikzpicture}
&\begin{tikzpicture}[x=4.4mm,y=8mm,baseline=(current bounding box.center)]
  \path (0,1.35);
  \draw[rounded corners,black!40] (1.55,-0.4) rectangle (3.45,1.3);
  \gpgraph{0,...,3}{1,2,3}{0,1,3}
  \gpcols{0,...,3}
  \node[font=\small] at (1.5,-1.4) {$*$};
\end{tikzpicture}
&\begin{tikzpicture}[x=4.4mm,y=8mm,baseline=(current bounding box.center)]
  \path (0,1.35);
  \draw[rounded corners,black!40] (2.55,-0.4) rectangle (4.45,1.3);
  \gpgraph{0,...,4}{1,2,3,4}{0,1,2,4}
  \gpcols{0,...,4}
  \node[font=\small] at (2.0,-1.4) {$S^{2}$};
\end{tikzpicture}
&\begin{tikzpicture}[x=4.4mm,y=8mm,baseline=(current bounding box.center)]
  \path (0,1.35);
  \draw[rounded corners,black!40] (3.55,-0.4) rectangle (5.45,1.3);
  \gpgraph{0,...,5}{1,2,3,4,5}{0,1,2,3,5}
  \gpcols{0,...,5}
  \node[font=\small] at (2.5,-1.4) {$*$};
\end{tikzpicture}
&\begin{tikzpicture}[x=4.4mm,y=8mm,baseline=(current bounding box.center)]
  \path (0,1.35);
  \draw[rounded corners,black!40] (4.55,-0.4) rectangle (6.45,1.3);
  \gpgraph{0,...,6}{1,2,3,4,5,6}{0,1,2,3,4,6}
  \gpcols{0,...,6}
  \node[font=\small] at (3.0,-1.4) {$S^{3}$};
\end{tikzpicture}\\[8pt]
$\Gamma_4$
&\begin{tikzpicture}[x=4.4mm,y=8mm,baseline=(current bounding box.center)]
  \path (0,1.35);
  \draw[rounded corners,black!40] (0.55,-0.4) rectangle (2.45,1.3);
  \gpgraph{0,...,2}{1,2}{0}
  \gpcols{0,...,2}
  \node[font=\small] at (1.0,-1.4) {$*$};
\end{tikzpicture}
&\begin{tikzpicture}[x=4.4mm,y=8mm,baseline=(current bounding box.center)]
  \path (0,1.35);
  \draw[rounded corners,black!40] (1.55,-0.4) rectangle (3.45,1.3);
  \gpgraph{0,...,3}{1,2,3}{0,1}
  \gpcols{0,...,3}
  \node[font=\small] at (1.5,-1.4) {$*$};
\end{tikzpicture}
&\begin{tikzpicture}[x=4.4mm,y=8mm,baseline=(current bounding box.center)]
  \path (0,1.35);
  \draw[rounded corners,black!40] (2.55,-0.4) rectangle (4.45,1.3);
  \gpgraph{0,...,4}{1,2,3,4}{0,1,2}
  \gpcols{0,...,4}
  \node[font=\small] at (2.0,-1.4) {$S^{2}$};
\end{tikzpicture}
&\begin{tikzpicture}[x=4.4mm,y=8mm,baseline=(current bounding box.center)]
  \path (0,1.35);
  \draw[rounded corners,black!40] (3.55,-0.4) rectangle (5.45,1.3);
  \gpgraph{0,...,5}{1,2,3,4,5}{0,1,2,3}
  \gpcols{0,...,5}
  \node[font=\small] at (2.5,-1.4) {$S^{2}$};
\end{tikzpicture}
&\begin{tikzpicture}[x=4.4mm,y=8mm,baseline=(current bounding box.center)]
  \path (0,1.35);
  \draw[rounded corners,black!40] (4.55,-0.4) rectangle (6.45,1.3);
  \gpgraph{0,...,6}{1,2,3,4,5,6}{0,1,2,3,4}
  \gpcols{0,...,6}
  \node[font=\small] at (3.0,-1.4) {$*$};
\end{tikzpicture}\\[8pt]
\end{tabular}
\caption{The graphs $R_m(\Gamma_j)$ for $0\le m\le4$. 
Under each graph is the homotopy type of the independence complex.}\label{fig:gp-base}
\end{figure}
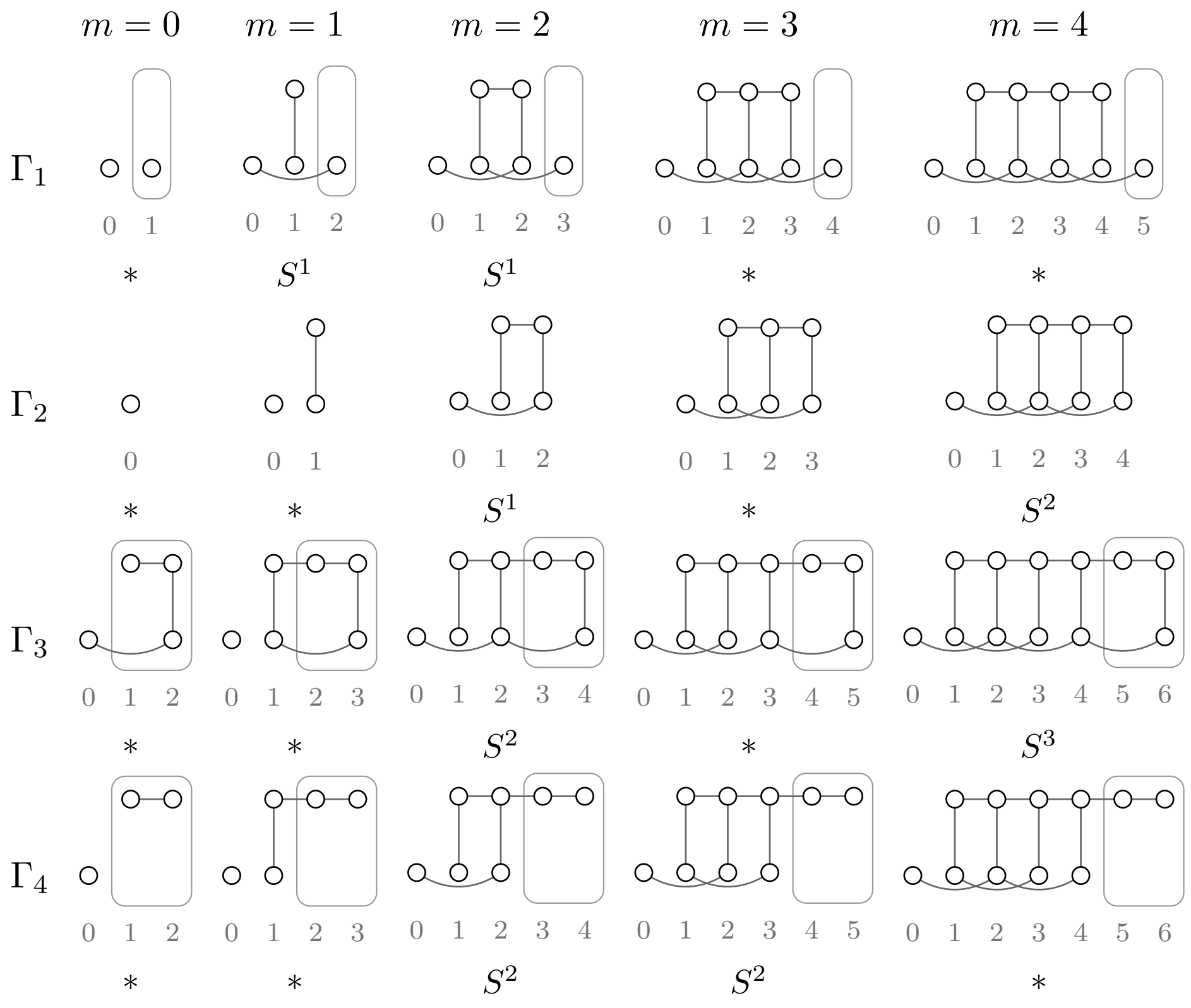

\begin{lemma}\label{lem:branches}
Let $n\ge5$, and write $n=5q+r$ with $0\le r\le 4$.
Then the homotopy types of $\Ind(B_1),\dots,\Ind(B_4)$ are as follows, where $*$ means
contractible.
\[
  \renewcommand{\arraystretch}{1.2}
  \begin{array}{@{}c*{5}{w{c}{3.6em}}@{}}
    \toprule
    r         & 0        & 1        & 2        & 3      & 4\\
    \midrule
    \Ind(B_1) & S^{3q-1} & S^{3q-1} & *        & *      & *\\
    \Ind(B_2) & S^{3q-2} & *        & S^{3q-1} & *      & *\\
    \Ind(B_3) & S^{3q-2} & *        & *        & S^{3q} & *\\
    \Ind(B_4) & S^{3q-3} & *        & *        & *      & S^{3q}\\
    \bottomrule
  \end{array}
\]
\end{lemma}
\begin{proof}
We write $B_j(n)$ for the graph $B_j$ in $\GP(n,2)$.
We apply \cref{lem:split} once to $B_1$, $B_3$ and $B_4$.
Then each $B_j$ becomes a graph $R_m(\Gamma_j)$, as can be checked in \cref{fig:cofibre}:
\begin{itemize}
\item $B_1=G\setminus\{v_0,u_1\}$. Applying \cref{lem:split} to $u_{n-1}$ leaves
  $R_{n-4}(\Gamma_1)$, where $\Gamma_1$ is the single vertex $b=v_{n-2}$.
\item $B_2=G\setminus\{u_0,u_1,u_2,v_0,v_1\}$ is $R_{n-3}(\Gamma_2)$, where $\Gamma_2$ is empty.
\item $B_3=G\setminus\{u_0,v_0,v_1,v_2,v_{n-2}\}$. Applying \cref{lem:split} to $u_2$ leaves
  $R_{n-6}(\Gamma_3)$, where $\Gamma_3$ is the path $a\,u_{n-1}\,c$ with $a=u_{n-2}$ and
  $c=v_{n-1}$.
\item $B_4=G\setminus\{u_0,u_1,v_0,v_1,v_2,v_3,v_{n-2},v_{n-1}\}$. Applying \cref{lem:split} to
  $u_3$ leaves $R_{n-7}(\Gamma_4)$, where $\Gamma_4$ is the edge $a\,u_{n-1}$ with $a=u_{n-2}$.
\end{itemize}
The vertices among $a,b,c$ that are not named are absent.
Hence
\begin{align*}
  \Ind(B_1(n))&\simeq\Sigma\Ind(R_{n-4}(\Gamma_1)) && (n\ge5),\\
  \Ind(B_2(n))&\simeq\Ind(R_{n-3}(\Gamma_2)) && (n\ge5),\\
  \Ind(B_3(n))&\simeq\Sigma\Ind(R_{n-6}(\Gamma_3)) && (n\ge6),\\
  \Ind(B_4(n))&\simeq\Sigma\Ind(R_{n-7}(\Gamma_4)) && (n\ge7).
\end{align*}
By \cref{prop:block}, $\Ind(R_m(\Gamma_j))\simeq\Sigma^3\Ind(R_{m-5}(\Gamma_j))$ for $m\ge5$, so
it remains to find $\Ind(R_m(\Gamma_j))$ for $0\le m\le4$.
These twenty graphs are drawn in \cref{fig:gp-base}.
Applying the steps (1)--(3) of the proof of \cref{prop:block} to each of them gives the
homotopy types in \cref{fig:gp-base}.
So by \cref{prop:block}, $\Ind(R_m(\Gamma_j))$ is as follows, where $t\ge0$.
\[
  \renewcommand{\arraystretch}{1.2}
  \begin{array}{@{}c*{5}{w{c}{3.6em}}@{}}
    \toprule
    m         & 5t   & 5t+1     & 5t+2     & 5t+3     & 5t+4\\
    \midrule
    \Gamma_1  & *    & S^{3t+1} & S^{3t+1} & *        & *\\
    \Gamma_2  & *    & *        & S^{3t+1} & *        & S^{3t+2}\\
    \Gamma_3  & *    & *        & S^{3t+2} & *        & S^{3t+3}\\
    \Gamma_4  & *    & *        & S^{3t+2} & S^{3t+2} & *\\
    \bottomrule
  \end{array}
\]
Substituting $m=n-4,n-3,n-6,n-7$ into the formulas for $\Ind(B_j(n))$ above gives the lemma.
The three graphs not covered by the formulas are the paths
\[
  B_3(5)=u_1u_2u_3u_4v_4,\qquad B_4(5)=u_2u_3u_4,\qquad B_4(6)=u_2u_3u_4u_5
\]
(see \cref{fig:small}).
They have homotopy types $\Ind(B_3(5))\simeq S^1$, $\Ind(B_4(5))\simeq S^0$, and $\Ind(B_4(6))\simeq *$, as
claimed.
\end{proof}

\begin{figure}[ht]
\centering
\begin{tikzpicture}[baseline=(current bounding box.north)]
  \gpcircle{5}{1,2,3,4}{4}
  \node[font=\small] at (0,-1.95) {$B_3(5)$};
\end{tikzpicture}
\qquad
\begin{tikzpicture}[baseline=(current bounding box.north)]
  \gpcircle{5}{2,3,4}{}
  \node[font=\small] at (0,-1.95) {$B_4(5)$};
\end{tikzpicture}
\qquad
\begin{tikzpicture}[baseline=(current bounding box.north)]
  \gpcircle{6}{2,3,4,5}{}
  \node[font=\small] at (0,-1.95) {$B_4(6)$};
\end{tikzpicture}
\caption{The graphs $B_3(5)$ and $B_4(5)$ in $\GP(5,2)$, and $B_4(6)$ in $\GP(6,2)$.}\label{fig:small}
\end{figure}
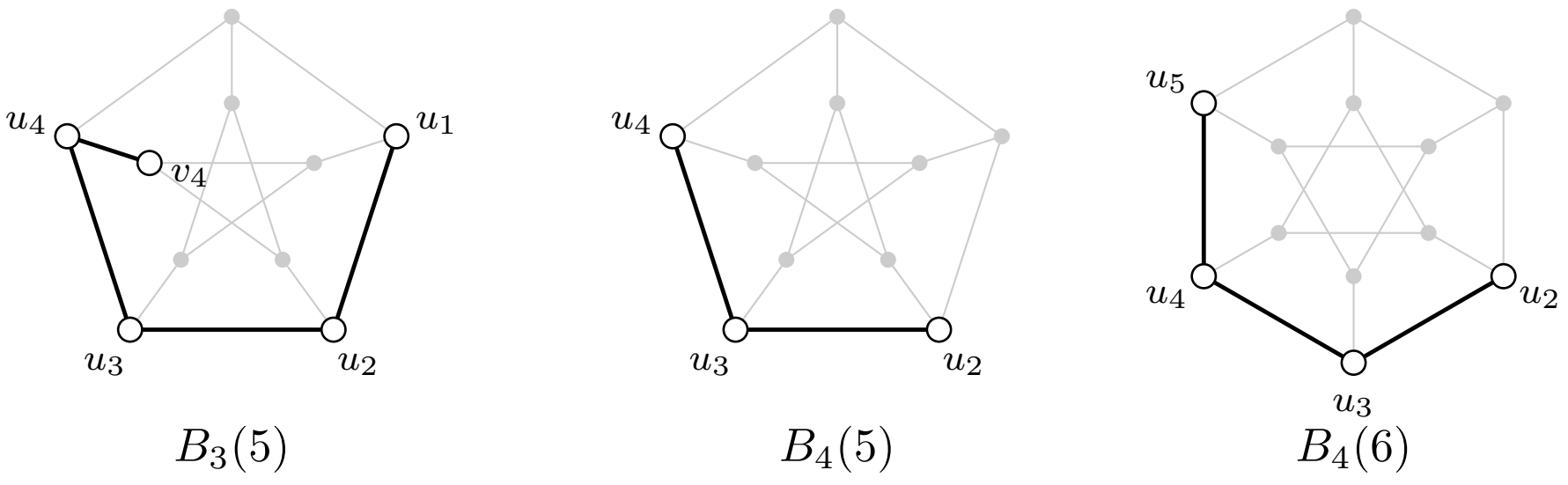

We can now determine $\Ind(\GP(n,2))$.

\begin{theorem}\label{thm:gp}
Let $n\ge5$, and write $n=5q+r$ with $r\in\{0,1,2,3,4\}$. Then
\[
  \Ind(\GP(n,2))\simeq
  \begin{cases}
    \bigvee_4S^{3q-1} & \text{if } r=0,\\
    S^{3q-2+r} & \text{if } r=1,2,3,4.
  \end{cases}
\]
Equivalently, $\Ind(\GP(n+5,2))\simeq\Sigma^3\Ind(\GP(n,2))$ for all $n\ge5$, where
$\Ind(\GP(5,2))\simeq\bigvee_4S^2$ and $\Ind(\GP(n,2))\simeq S^{n-4}$ for $n=6,7,8,9$.
\end{theorem}

\begin{proof}
The three cofibre sequences of \cref{fig:cofibre}, at $u_1$ in $G\setminus v_0$, at $v_1$ in
$G\setminus N[v_0]$, and at $v_0$ in $G$, are
\begin{align*}
  \Ind(B_2)&\hookrightarrow\Ind(B_1)\hookrightarrow\Ind(G\setminus v_0),\\
  \Ind(B_4)&\hookrightarrow\Ind(B_3)\hookrightarrow\Ind(G\setminus N[v_0]),\\
  \Ind(G\setminus N[v_0])&\hookrightarrow\Ind(G\setminus v_0)\hookrightarrow\Ind(G).
\end{align*}
We perform a case-by-case analysis applying \cref{lem:cofibre} and using the values in \cref{lem:branches}.

Let $r=1$ or $r=2$.
Then $\Ind(B_3)\simeq\Ind(B_4)\simeq*$, so $\Ind(G\setminus N[v_0])\simeq*$, and thus
$\Ind(G)\simeq\Ind(G\setminus v_0)$.
Since $\Ind(B_2)\simeq*$ for $r=1$ and $\Ind(B_1)\simeq*$ for $r=2$, we get
\[
  \Ind(G)\simeq\Ind(G\setminus v_0)\simeq
  \begin{cases}
    \Ind(B_1)\simeq S^{3q-1} & \text{if } r=1,\\
    \Sigma\Ind(B_2)\simeq S^{3q} & \text{if } r=2.
  \end{cases}
\]

Let $r=3$ or $r=4$.
Then $\Ind(B_1)\simeq\Ind(B_2)\simeq*$, so $\Ind(G\setminus v_0)\simeq*$, and thus
$\Ind(G)\simeq\Sigma\Ind(G\setminus N[v_0])$.
Since $\Ind(B_4)\simeq*$ for $r=3$ and $\Ind(B_3)\simeq*$ for $r=4$, we get
\[
  \Ind(G)\simeq\Sigma\Ind(G\setminus N[v_0])\simeq
  \begin{cases}
    \Sigma\Ind(B_3)\simeq S^{3q+1} & \text{if } r=3,\\
    \Sigma^2\Ind(B_4)\simeq S^{3q+2} & \text{if } r=4.
  \end{cases}
\]
Therefore $\Ind(G)\simeq S^{3q-2+r}$ for $r=1,2,3,4$.

If $r=0$, then $q\ge1$, and $\Ind(B_1),\dots,\Ind(B_4)$ are spheres of dimensions $3q-1$,
$3q-2$, $3q-2$, $3q-3$.
In the first two cofibre sequences, the first complex is a sphere of one dimension lower than
the second, so by \cref{lem:cofibre}(c),
\[
  \begin{aligned}
    \Ind(G\setminus v_0)&\simeq S^{3q-1}\vee\Sigma S^{3q-2}=\bigvee_2S^{3q-1},\\
    \Ind(G\setminus N[v_0])&\simeq S^{3q-2}\vee\Sigma S^{3q-3}=\bigvee_2S^{3q-2}.
  \end{aligned}
\]
The same dimension gap holds in the third sequence, hence
\[
  \Ind(G)\simeq\bigvee_2S^{3q-1}\vee\Sigma\bigvee_2S^{3q-2}=\bigvee_4S^{3q-1}.
  \qedhere
\]
\end{proof}

The Petersen graph $\GP(5,2)$ is the Kneser graph of $2$-subsets of $\{1,\dots,5\}$, and its
independence complex was determined by Barmak \cite[Thm.~4.11]{Barmak13}, as part of the
family of Kneser graphs on $2$-subsets.

\section{The perfect matching complex of \texorpdfstring{$\GP(n,2)$}{GP(n,2)}}\label{sec:gp-pm}

By \cref{thm:gp}, the independence complexes of the graphs $\GP(n,2)$ are spheres or wedges of
four spheres, and they are periodic in $n$.
We end with a discussion of the perfect matching complexes of these graphs.
For the Petersen graph we determine $\Mp(\GP(5,2))$ by hand (\cref{prop:petersen}).
For $6\le n\le14$ we computed the homology of $\Mp(\GP(n,2))$, which shows no evident pattern
(\cref{tab:gp-pm}).

\begin{proposition}\label{prop:petersen}
The Petersen graph $\GP(5,2)$ has six perfect matchings, and any two of them share exactly one
edge. Moreover, $\Mp(\GP(5,2))\simeq\bigvee_{10}S^1$.
\end{proposition}

\begin{proof}
It is easy to see that the perfect matchings are $M=\{u_0v_0,\dots,u_4v_4\}$ and,
for $i\in\Z_5$,
\[
  P_i=\{u_iv_i,\ u_{i+1}u_{i+2},\ u_{i+3}u_{i+4},\ v_{i+2}v_{i+4},\ v_{i+1}v_{i+3}\},
\]
drawn below.
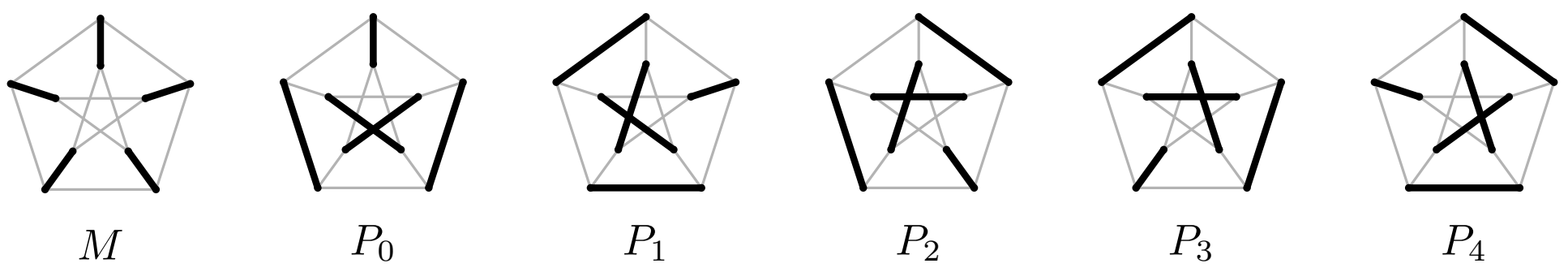
\begin{center}
\petersenpm{u0/v0,u1/v1,u2/v2,u3/v3,u4/v4}{$M$}\hfill
\petersenpm{u0/v0,u1/u2,u3/u4,v2/v4,v1/v3}{$P_0$}\hfill
\petersenpm{u1/v1,u2/u3,u4/u0,v3/v0,v2/v4}{$P_1$}\hfill
\petersenpm{u2/v2,u3/u4,u0/u1,v4/v1,v3/v0}{$P_2$}\hfill
\petersenpm{u3/v3,u4/u0,u1/u2,v0/v2,v4/v1}{$P_3$}\hfill
\petersenpm{u4/v4,u0/u1,u2/u3,v1/v3,v0/v2}{$P_4$}
\end{center}
We have $M\cap P_i=\{u_iv_i\}$, $P_i\cap P_{i+1}=\{v_{i+2}v_{i+4}\}$ and
$P_i\cap P_{i+2}=\{u_{i+3}u_{i+4}\}$, so any two of them share exactly one edge.
These fifteen intersections are distinct, so no edge lies in three perfect matchings.
Hence $\Mp(\GP(5,2))$ is covered by six simplices $M,P_0,\dots,P_4$, any two of which meet in a $0$-simplex, and
no three of which meet.
By the nerve lemma \cite[Sect.~15.4]{KozlovBook}, $\Mp(\GP(5,2))$ is homotopy equivalent to the
nerve of this cover. The nerve is the complete graph $K_6$, which is connected and has Euler
characteristic $6-15=-9$, so it is homotopy equivalent to $\bigvee_{10}S^1$.
\end{proof}

For $6\le n\le14$, a computation with SageMath \cite{Sage} gives the reduced integral homology
of $\Mp(\GP(n,2))$ in \cref{tab:gp-pm}.
It shows no periodicity in $n$ and no obvious pattern in the ranks.
For $n=10$ it has torsion, so $\Mp(\GP(10,2))$ is not homotopy equivalent to a wedge of spheres.

\begin{table}[ht]
\centering
\renewcommand{\arraystretch}{1.2}
\begin{tabular}{@{}cl@{\qquad}cl@{}}
  \toprule
  $n$ & $\tilde H_*(\Mp(\GP(n,2));\Z)$ & $n$ & $\tilde H_*(\Mp(\GP(n,2));\Z)$\\
  \midrule
  $5$ & $\tilde H_1\cong\Z^{10}$ & $10$ & $\tilde H_2\cong\Z/2$, $\tilde H_3\cong\Z^{20}$\\
  $6$ & $0$ & $11$ & $\tilde H_3\cong\Z^{22}$\\
  $7$ & $\tilde H_2\cong\Z^7$ & $12$ & $\tilde H_3\cong\Z^{13}$, $\tilde H_4\cong\Z^3$\\
  $8$ & $\tilde H_2\cong\Z^{10}$ & $13$ & $\tilde H_4\cong\Z^{13}$\\
  $9$ & $\tilde H_2\cong\Z^9$ & $14$ & $\tilde H_4\cong\Z^{65}$, $\tilde H_5\cong\Z^2$\\
  \bottomrule
\end{tabular}
\caption{The nonzero reduced integral homology groups of $\Mp(\GP(n,2))$ for $5\le n\le14$.}
\label{tab:gp-pm}
\end{table}

\bibliographystyle{amsalpha}
\bibliography{antiprisms}

\end{document}